\documentclass[12pt]{amsart}

\usepackage[colorlinks=true, pdfstartview=FitV, linkcolor=blue,
citecolor=blue, urlcolor=blue]{hyperref}

\usepackage{mathtools}
\usepackage{amsfonts}
\usepackage{amsmath}
\usepackage{amssymb} 
\usepackage{amsthm}

\usepackage{times}
\usepackage[T1]{fontenc}
\usepackage{enumitem}
\usepackage{setspace}
\usepackage{microtype}
\usepackage{cite}

\allowdisplaybreaks

\usepackage[margin=1.2in]{geometry}

\newtheorem{theorem}{Theorem}[section]

\theoremstyle{definition}

\theoremstyle{remark}
\newtheorem{remark}[theorem]{Remark}

\numberwithin{equation}{section}

\newcommand{\R}{\mathbb R}
\newcommand{\N}{\mathbb N}

\newcommand{\Z}{\mathbb Z}

\newcommand{\rn}{{{\mathbb R}^n}}
\newcommand{\rd}{{{\mathbb R}^d}}

\DeclareMathOperator{\supp}{supp}

\DeclareMathOperator*{\esssup}{ess\,sup}

\newcommand{\D}{\mathcal D}
\newcommand{\Ss}{\mathcal S}
\newcommand{\B}{\mathcal{B}}
\newcommand{\Q}{\mathcal{Q}}

\newcommand{\Zz}{\mathcal{Z}}
\newcommand{\loc}{\text{loc}}

\newcommand{\F}{\mathcal{F}}
\newcommand{\Rdf}{\mathcal{R}}

\DeclareMathOperator{\op}{op}

\newcommand{\A}{\mathcal A}
\newcommand{\K}{\mathcal K}
\newcommand{\W}{\mathcal W}

\DeclareMathOperator{\clconv}{\overline{conv}}

\newcommand{\vecf}{\mathbf{f}}
\newcommand{\vecg}{\mathbf{g}}

\newcommand{\Sd}{\mathcal{S}_d}

\newcommand{\pp}{{p(\cdot)}}

\def\Xint#1{\mathchoice
   {\XXint\displaystyle\textstyle{#1}}%
   {\XXint\textstyle\scriptstyle{#1}}%
   {\XXint\scriptstyle\scriptscriptstyle{#1}}%
   {\XXint\scriptscriptstyle\scriptscriptstyle{#1}}%
   \!\int}
\def\XXint#1#2#3{{\setbox0=\hbox{$#1{#2#3}{\int}$}
     \vcenter{\hbox{$#2#3$}}\kern-.5\wd0}}

\def\avgint{\Xint-}

\begin{document}

\title{Recent results on matrix weighted  norm inequalities}

\author{David Cruz-Uribe, OFS}
\address{Dept. of Mathematics \\
University of Alabama \\
 Tuscaloosa, AL 35487, USA}
\email{dcruzuribe@ua.edu}


\date{August 18, 2025}

\thanks{The author is partially supported by a Simons Foundation
  Travel Support for Mathematicians  Grant and by NSF Grant
  DMS-2349550.  Earlier versions of this article were based on a talk
  given at the Conference on Research on Real, Complex and Functional Analysis
  in Kyoto, Japan, October 26--27, 2022; see arXiv:2304.03887. The author would like to thank
  Yoshihiro Sawano for his invitation and generous support to speak at
  this conference.  This substantially revised version is based on
  talks  given at the conference Expanding Pathways in Harmonic
  Analysis, Washington University in St. Louis, April 5--6, 2025, and
  the conference on Advances in Nonlinear Analysis, held
  at Rzeszow University of Technology, Poland, June 30--July 2,
  2025. The author is grateful for the opportunity to speak at this
  conference held in honor of the 75th birthday of Professor J\'ozef
  Bana\'s: {\em Sto Lat!}}

\begin{abstract}
In this paper we give an overview of recent work on matrix weights,
with particular emphasis on convex body sparse domination for singular
integrals, Rubio de Francia extrapolation, and Jones factorization.
To provide context and  motivation, we survey the comparable results
in the scalar weighted case.  
\end{abstract}

\maketitle

\section{Introduction}
The goal of this article is to provide an overview of recent work on
matrix weights and matrix weighted norm inequalities.  In the scalar
case, 
one weight norm inequalities have been extensively studied since the
work of Muckenhoupt and others in the
1970s and there is now a rich and well developed theory.   (See~\cite{duoandikoetxea01,garcia-cuerva-rubiodefrancia85,grafakos08b}.)  Central to their study is the
Muckenhoupt $A_p$ condition.  It gives a necessary and sufficient
condition for weighted strong and weak type inequalities for many of the
operators of classical harmonic analysis:  maximal operators, singular
integral operators, square functions.
 They also
play an important role in other fields, such as the study of
degenerate elliptic PDEs.  (See, for
instance,~\cite{MR643158,MR847996}.) 
A parallel theory of
``off-diagonal'' inequalities, primarily for fractional integrals
(i.e., Riesz potentials), has been developed using the
$A_{p,q}$ weights of Muckenhoupt and
Wheeden~\cite{muckenhoupt-wheeden74,CruzUribe:2016ji}.

The Muckenhoupt $A_p$ have a rich structure.  One of the deepest results in
this area is the Jones factorization
theorem~\cite{Jones}, which gives a complete characterization of the
class $A_p$ in terms of the simpler class $A_1$. 
Closely related to factorization is the Rubio de
Francia extrapolation theorem~\cite{rubiodefrancia84,MR2797562}, which, in its
simplest form, shows that if an operator is bounded on weighted $L^2$
for weights in $A_2$, then
it is bounded on weighted $L^p$ for all $p$, $1<p<\infty$, for weights
in $A_p$.   Rubio de
Francia's colleague Antonio Cordoba~\cite{garcia-cuerva87} summarized
this remarkable result by
saying that it showed that  {\em there are no $L^p$ spaces, only
  weighted $L^2$.}  

A more recent problem has been to  determine the sharp
dependence in weighted norm inequalities on the $A_p$ characteristic
$[w]_{A_p}$, the constant that comes from the Muckenhoupt condition.
This was question was first settled for the maximal operator by
Buckley~\cite{buckley93}.  In a major advance, 
Hyt\"onen~\cite{MR2912709} proved the sharp dependence for singular
integrals is $[w]_{A_p}^{\max\{1, \frac{1}{p-1}\}}$.   Central to his proof was a quantitative version of Rubio de
Francia extrapolation, which allowed him to reduce the problem to
proving the case $p=2$: that is, that the sharp dependence is $[w]_{A_2}$.  For this reason, the
problem was referred to as the ``$A_2$ conjecture.''

Since the 1990s, there has been a great deal of interest in extending
the theory of scalar weights to the setting of matrix weights:  that
is, $d\times d$ measurable matrix functions $W$ that are self-adjoint,
positive semidefinite, and act on vector-valued functions $\vecf$.
This problem was first studied by Nazarov, Treil and
Volberg~\cite{MR1428988,MR1428818,MR1423034,MR1478786}, who asked if
the strong type
weighted norm inequalities for the Hilbert transform on the real line could be extended to
the matrix setting.  Beyond
the intrinsic interest of this problem, their original motivation came
from problems in the study of multivariate random stationary
processes and  Toeplitz operators acting on
vector-valued Hardy spaces.

They defined a class of matrix weights, again denoted by $A_p$, and
showed that the Hilbert transform is bounded with respect to this
class.  This result was later extended to general singular integral
operators by Christ and Goldberg~\cite{MR2015733,MR1813604}. More
recently, other operators have been considered, such as square
functions~\cite{MR3936542,MR4159390},  and
commutators~\cite{MR4454483}.  Matrix $A_{p,q}$ weights were also
defined and studied to prove bounds for fractional
integrals~\cite{MR4030471}.  Originally, the focus was on strong type
inequalities, but more recently weak type inequalities have been defined
and proved~\cite{MR4269407,MR4859193}.  

After the proof of the $A_2$ conjecture for scalar weights, it was conjectured that the
same result held for matrix weights.   The
best known result is that the constant is bounded above by
$[W]_{\A_p}^{1+\frac{1}{p-1}-\frac{1}{p}}$, which, when $p=2$ gives an
exponent of $\frac{3}{2}$ rather than the conjectured value of $1$.
The case $p=2$ was proved by Nazarov, Petermichl, Treil and
Volberg~\cite{MR3689742}, and for all $p$ in~\cite{MR3803292}.  
Very recently, Domelevo, Petermichl, Treil and
Volberg~\cite{DPTV-2024} show that this conjecture is false and that
  when $p=2$, $\frac{3}{2}$ is the best possible exponent.  The
  question of sharp exponents when $p\neq 2$ is still open.
  
The structure of matrix $A_p$ weights is more complicated than that of
scalar weights.  Bownik~\cite{Bow} proved that some properties of scalar
weights, such as left openness (that is, if $w\in A_p$ for some $p>1$,
there exists $\epsilon>0$ such that $w\in A_{p-\epsilon}$) no longer
hold in the matrix case.  In the 1990s
Nazarov and Treil~\cite{MR1428988} posed two fundamental problems:
extend the Jones factorization theorem and the Rubio de Francia
extrapolation theorem to matrix weights.  They wrote:
\begin{quote}
    {\em Actually, the whole theory of scalar $(A_p)$-weights  can be
      transferred to the matrix case except two results...the Peter
      Jones factorization and the Rubio-de-Francia extrapolation
      theory.  But today (May 1, 1996) we do not know what the
      analogues of these two things are in high dimensions.}
    \end{quote}
These problems were finally
solved in 2022~\cite{mb-2022}.  Initially, one motivation for proving
extrapolation for matrix weights was to use it to prove the matrix
$A_2$ conjecture.  However, applying extrapolation to the sharp bound
for $p=2$ yields a worse estimate than is known for other values of
$p$.   It is an interesting problem to understand why extrapolation
fails in the matrix case to give the best possible dependence.

\smallskip

In the remainder of this article we will give a more detailed overview
of the theory of matrix weights, emphasizing the similarities and differences with
the scalar case.  Though we will provide detailed references, our goal
is not to provide a comprehensive history of the field; we will instead focus on
some specific results which we believe are at the heart of recent
research, particularly with the past five or ten years.  

It is organized
as follows.  In Section~\ref{section:scalar} we provide an overview of
scalar weighted norm inequalities, focusing on those results which are
essential for understanding the subsequent development of the matrix
theory.  As a first step we reformulate strong type inequalities using
the approach of treating the weight as a ``multiplier'' instead of as
a ``measure.''  This makes the transition to matrix weights easier in
later sections, but we caution the reader to be mindful of this when
comparing the results given here with the original literature.  We introduce both the classical weak type inequality, and
the so-called multiplier weak type inequalities, which were the
motivation for the definition of matrix weighted weak type
inequalities. We also include some brief remarks about fractional
integrals and fractional maximal operators.

In Section~\ref{section:jones-rubio} we review the scalar theory of
extrapolation and factorization.  We include some generalizations,
including to $A_p$ weights defined with respect 
to bases other than cubes, and off-diagonal extrapolation.

We then turn to the theory of matrix weights.  
In Section~\ref{section:matrix} we consider matrix weighted
inequalities for singular integrals.  We first discuss the three
equivalent definitions of matrix $A_p$.  We then give the strong
type inequalities of Christ and
Goldberg~\cite{MR2015733,MR1813604} and the results
proved by Nazarov,
Petermichl, Treil, and Volberg~\cite{MR3689742}  using the technique of
convex body sparse domination.  We will then consider weak type
inequalities.   We will show how the fine properties of
matrix $A_p$ weights can be combined with weak type inequalities to
prove strong type inequalities by interpolation.     We will also briefly
discuss some result for fractional integrals.  

In Section~\ref{section:convex} we will explore the techniques
underlying the proofs of Rubio de Francia extrapolation and Jones
factorization in the matrix setting.  These proofs rely on the theory
of harmonic analysis on convex set-valued functions which allows us to
generalize the proofs in the scalar case to the matrix case.  We will
then discuss some recent generalizations of extrapolation in the
matrix setting.  As an application of the theory of convex set-valued
functions, we will sketch a new and extremely elementary proof of
strong type bounds for singular integrals.

Finally, in Section~\ref{section:beyond} we
will briefly consider some
other results involving matrix weights.  This is not intended to be a
comprehensive summary.  Rather, we want to highlight some recent work
that we think is of interest, and give references for those wanting 
to explore further.

Throughout this paper we will use the following notation.  In
Euclidean space the constant $n$ will denote the dimension of $\R^n$,
which will be the domain of our functions.  The value $d$ will denote
the dimension of vector and set-valued functions.  For
$1\leq p \leq \infty$, $L^p(\R^n)$ will denote the Lebesgue space of
scalar functions, and $L^p(\R^n,\rd )$ will denote the Lebesgue space
of vector-valued functions.

By a cube we will always mean a cube in $\rn$ with sides parallel to
the coordinate axes.  The Lebesgue measure of a cube, or of any
arbitrary set $E$, will be denoted
by $|E|$.  By a weight we mean a nonnegative, locally
integrable function that is positive except on a set of measure $0$.
We define $w(E) = \int_E w(x)\,dx$, and we let $\avgint_Q w(x)\,dx =
|Q|^{-1}w(Q)$.  

Given $v=(v_1,\ldots,v_d)^t \in \rd $, the Euclidean norm of $v$ will
be denoted by $|v|$; from context there should be no confusion with
the notation for the Lebesgue measure of a set.  The closed unit ball
$\{ v \in \rd : |v| \leq 1\}$ will be denoted by ${\mathbf B}$.
Matrices will be $d\times d$ matrices with real-valued entries.  The
set of all $d\times d$, symmetric, positive semidefinite matrices will
be denoted by $\Sd$.

The letters $C$ and $c$ will denote constants whose value may change
at each appearance which will depend implicitly on the parameters
involved.  If we want to show explicit dependence, we will write, for
instance, $C(n,p)$.  Given two quantities $A$ and $B$, we will write
$A \lesssim B$, or $B\gtrsim A$ if there is a constant $c>0$ such that
$A\leq cB$.  If $A\lesssim B$ and $B\lesssim A$, we will write
$A\approx B$.

 \section{Scalar weighted norm inequalities}
 \label{section:scalar}

 To understand the motivation for the study of matrix
 weights, we first review some basic results in the theory of scalar
 weights.   We will consider the weights that
 satisfy the Muckenhoupt $A_p$ condition, introduced
 in~\cite{muckenhoupt72}.  Given $1<p<\infty$, $w\in A_p$ if
\[ [w]_{A_p} = \sup_Q \left(\avgint_Q w(x)\,dx\right) \left(\avgint_Q
  w(x)^{1-p'}\,dx\right)^{p-1} < \infty, \]
where the supremum is taken over all cubes.  A weight $w$ is in $A_1$ if
\[ [w]_{A_1} = \sup_{Q} \esssup_{Q\ni x}  w(x)^{-1} \avgint_Q
w(y)\,dy < \infty. \]
The quantity $[w]_{A_p}$ is referred to as the $A_p$ characteristic of
a weight. 

Here, we are going to depart from the classical definition which
treats the weight $w$ as a measure,  that is, the inequalities are
with respect to the measure $wdx$.   Instead,
we will consider the weight as a multiplier:  hereafter, for $1\leq
p\leq \infty$, and a nonnegative, measurable function $w$, we define
$L^p(w)$ to be the class of measurable functions such that $wf \in
L^p(\rn)$ with norm $\|f\|_{L^p(w)} = \|wf\|_{L^p(\rn)}$.  
  This approach to weights was introduced by Muckenhoupt
and Wheeden~\cite{muckenhoupt-wheeden74} to study fractional integrals
(see below) but has also been used to extend the theory of
weights to other function spaces, such as the variable Lebesgue spaces.
(See, for instance,~\cite{MR2927495}.)  

We define the class of weights $\A_p$ as follows.  If $1<p<\infty$,
then $w\in \A_p$ if
\[ [w]_{\A_p} =\sup_Q \bigg(\avgint_Q w(y)^p\,dy\bigg)^{\frac{1}{p}}
  \bigg(\avgint_Q w(y)^{-p'}\,dy\bigg)^{\frac{1}{p'}} < \infty;\]
if $p=1$, 
\[ [w]_{\A_1} =\sup_Q \bigg(\avgint_Q w(y)\,dy\bigg) \esssup_{x\in Q}  w(x)^{-1}
  < \infty; \]
and if $p=\infty$,
\[ [w]_{\A_\infty} =\esssup_{x\in Q}  w(x)\sup_Q \bigg(\avgint_Q w(y)^{-1}\,dy\bigg) 
  < \infty. \]

\begin{remark}
The two conditions $A_1$ and $\A_1$ are the same.  For $1<p<\infty$,
the conditions $A_p$ and $\A_p$ are equivalent  under the mapping
$w\mapsto w^p$:  $w\in \A_p$ if and only if $w^p \in A_p$, and
$[w]_{\A_p}^p=[w^p]_{A_p}$.   When $p=\infty$ we have that $w\in \A_\infty$
  if and only if $w^{-1}\in \A_1=A_1$.  Note that this definition of
  $\A_\infty$  is quite different from the classical $A_\infty$
  condition.
\end{remark}

\begin{remark}
  All of these definitions can be written in a single combined form
  using norms instead of integrals:  for $1\leq p \leq \infty$, $w\in
  \A_p$ if and only if
  \[ [w]_{\A_p} = \sup_Q |Q|^{-1}
    \|w\chi_Q\|_{L^p}\|w^{-1}\chi_Q\|_{L^{p'}}<\infty.  \]
  This formulation also has the advantage that it extends naturally to
  other Banach functions spaces:  see~\cite{MR2927495,
    nieraeth2025muckenhouptcondition, MR4632742}.
\end{remark}

The Muckenhoupt $\A_p$ weights arise naturally in the study of the Hardy-Littlewood maximal
operator,
\[ Mf(x) = \sup_Q \avgint_Q |f(y)|\,dy \cdot \chi_Q(x). \]
The following theorem is
due to Muckenhoupt~\cite{muckenhoupt72}.  For a more elementary proof
of the strong type inequality due to Christ and Fefferman, see~\cite{MR684636,dcu-paseky}.

\begin{theorem} \label{thm:scalar-max}
Fix  $1\leq p \leq \infty$.  Given a weight $w$, the following are equivalent:

  \begin{enumerate}

\item $w\in \A_p$;

\item the maximal operator satisfies the weak  $(p,p)$ inequality
  \[ \sup_{\lambda>0} \lambda \| \chi_{ \{ x : Mf(x)>\lambda\}}
      \|_{L^p(w)}
      \leq
      C\|f\|_{L^p(w)}; \]

  \item if additionally $p>1$, then
    \[ \|Mf\|_{L^p(w)} \leq \|f\|_{L^p(w)}. \]
  \end{enumerate}
\end{theorem}

\begin{remark}
  If $p<\infty$, the weak $(p,p)$ inequality is usually written
\[ \sup_{\lambda>0} \lambda^p w^p(\{ x\in \rn : Mf(x) > \lambda \})
  \leq
C\int_\rn |w(x)f(x)|^p \,dx, \]
and when  $1<p<\infty$ the strong
$(p,p)$ inequality is written
\[ \int_\rn (w(x)Mf(x))^p \,dx
  \leq
  C \int_\rn |w(x)f(x)|^p \,dx. \]
One advantage of the $\A_p$ condition is that we can naturally extend
the norm inequalities to  the case $p=\infty$.  For the maximal
operator this was done in
Muckenhoupt's original work, but seems to have been mostly
overlooked since then. 
\end{remark}

\smallskip

Muckenhoupt $\A_p$ weights also control weighted norm inequalities for
singular integrals.    A Calder\'on-Zygmund singular integral $T$
is a bounded
operator on $L^2(\rn)$ for which there exists a kernel $K(x,y)$, defined on
$\rn\times \rn\setminus \Delta$, where $\Delta =\{ (x,x) : x \in \rn\}$, such that
if $f \in L_c^\infty(\rn)$ and $x\not\in \supp(f)$, then
\[ Tf(x) = \int_\rn K(x,y)f(y)\,dy.  \]
The kernel satisfies the size and regularity conditions
\[ |K(x,y)| \leq \frac{C}{|x-y|^n}, \]
\[ |K(x+h,y)-K(x,y)|+|K(x,y+h)-K(x,y)|
  \leq
  C\frac{|h|^\delta}{|x-y|^{n+\delta}}, \]
where $|x-y|>2|h|$.  We say that the kernel is nondegenerate if there
exist $a>0$ and a unit vector $v\in \rn$ such that if $x,\,y\in \rn$ satisfy
$x-y=tv$, $t\in \R$, then
\[ |K(x,y)| \geq \frac{a}{|x-y|^n}. \]
For example, this condition is satisfied by the Hilbert transform or
any of the Riesz transforms.

In the following result, the sufficiency is due to Hunt,
Muckenhoupt and Wheeden~\cite{hunt-muckenhoupt-wheeden73} for the
Hilbert transform and Coifman and Fefferman~\cite{coifman-fefferman74}
for general singular integrals; the necessity  is due to
Stein~\cite{stein93}.   See also~\cite{duoandikoetxea01,garcia-cuerva-rubiodefrancia85,
  grafakos08b}.

\begin{theorem} \label{thm:scalar-sio}
Given a Calder\'on-Zygmund singular integral operator $T$, if
$1 \leq p<\infty$ and  $w\in \A_p$, then the weak $(p,p)$ inequality 
\[ \sup_{\lambda>0} \lambda \|\chi_{\{ x : |Tf(x)| > \lambda \}}\|_{L^p(w)}
  \leq
  C  \|f\|_{L^p(w)}. \]
holds.  Additionally, if $1<p<\infty$, then the strong $(p,p)$
inequality
\[ \|Tf\|_{L^p(w)}
  \leq
  C \|f\|_{L^p(w)} \]
holds. If $T$ is nondegenerate, then the strong and weak type
inequalities imply that $w\in \A_p$.
\end{theorem}

The weak $(p,p)$ inequality for singular integrals can be proved using
kernel estimates and the good/bad decomposition of Calder\'on and
Zygmund. (See, for instance,~\cite{garcia-cuerva-rubiodefrancia85}.)
The strong type inequality was originally proved by comparing the norm
of the singular integral operator to that of the maximal operator.
Coifman and Fefferman~\cite{coifman-fefferman74} proved that given
$p$, $0<p<\infty$, and $w\in \A_q$ for any $q$, $1\leq q\leq \infty$, there
exists a constant depending on $[w]_{\A_q}$ such that
\begin{equation} \label{eqn:TM}
  \|T^*f\|_{L^p(w)} \leq C \|Mf\|_{L^p(w)}. 
\end{equation}
Here, $T^*$ is the maximal singular integral, defined by
\[ T^* f(x) = \sup_{\epsilon>0}|T_\epsilon f(x)|
  = \sup_{\epsilon>0}
  \bigg|\int_{|x-y|>\epsilon} K(x,y)f(y)\,dy\bigg|, \]
which dominates the singular integral pointwise.  They proved this
inequality by
using a good-$\lambda$ inequality:  they showed that there exists
$\delta>0$ such that for every $\gamma,\,\lambda>0$ and for every cube $Q$,
\begin{equation} \label{eqn:good-lambda-sio}
 w(\{ x\in Q : T^*f(x) >2\lambda, Mf(x)\leq \gamma\lambda\})
  \leq
  C\gamma^\delta w(Q). 
\end{equation}
An alternative proof of \eqref{eqn:TM} using the sharp maximal
operator was given by
Journ\'e~\cite{MR706075}; see also Alvarez and
P\'erez~\cite{alvarez-perez94}.  

\smallskip

The constants in the norm inequalities in
Theorems~\ref{thm:scalar-max} and~\ref{thm:scalar-sio} depend on the
$\A_p$ characteristic, but the original proofs were qualitative and
did not give the quantitative dependence on it.  Since the 1990s there has been a
great deal of interest in determining the best constant in the strong
$(p,p)$ inequalities in terms of $[w]_{\A_p}$.  This question was
first considered by Buckley~\cite{buckley93}.  He proved
the sharp bound for the maximal operator:  for $1<p<\infty$,
\[  \|Mf\|_{L^p(w)} \leq C[w]_{\A_{p}}^{p'} \|f\|_{L^p(w)}. \]
He also  gave a lower bound for
singular integrals that was conjectured to be sharp.  This conjecture
became the subject of concerted effort when Astala, Iwaniec and
Saksman~\cite{MR1815249} proved that sharp regularity results for
solutions of the Beltrami equation hold provided that the
Beurling-Ahlfors operator satisfies
$\|Tf\|_{L^p(w)} \leq C [w]_{\A_p}^p\|f\|_{L^p(w)}$ for $p\geq 2$.  This
conjecture was extended to all Calder\'on-Zygmund operators:  for any
$p>1$,
\begin{equation} \label{eqn:sharp} 
\|Tf\|_{L^p(w)} \leq C(n,T,p) [w]_{\A_p}^{\max(p,p')}\|f\|_{L^p(w)}.
\end{equation}
By the Rubio de Francia extrapolation theorem, it suffices to prove
this conjecture for $p=2$.  Results for specific singular integrals
(e.g., the Beurling-Ahlfors operator, the Hilbert transform, Riesz
transforms) were obtained by Petermichl and
Volberg~\cite{MR2354322,petermichl08, petermichl-volberg02}.  Partial
results for the general problem were obtained by several authors (see
\cite{lacey-petermichl-reguera2010,dcu-martell-perez}); the $A_2$
conjecture was finally proved by
Hyt\"onen~\cite{MR2912709,hytonen-perez-treil-volbergP}.

  Sharp bounds for the weak type inequalities for maximal operators
  and singular integral operators are also known.  For the maximal
  operator, it follows immediately from the proof of the weak $(p,p)$
  inequality, $1\leq p < \infty$, that
  \[ \sup_{\lambda>0} \lambda \| \chi_{ \{ x : Mf(x)>\lambda\}}
      \|_{L^p(w)}
      \leq
      C[w]_{\A_p}\|f\|_{L^p(w)}; \]
    this bound was first proved by Muckenhoupt~\cite{muckenhoupt72} who also
    proved that it was the best possible.   The sharp bound for
    singular integral operators is
     \[ \sup_{\lambda>0} \lambda \| \chi_{ \{ x : |Tf(x)|>\lambda\}}
      \|_{L^p(w)}
      \leq
      C[w]_{\A_p}^p\|f\|_{L^p(w)}; \]
    this bound was proved by Hyt\"onen, {\em et al.}~\cite{MR2993026}.

  \smallskip

  Hyt\"onen's proof of the $A_2$ conjecture was very complicated.
  Soon after it was proved, however, 
Lerner and Nazarov~\cite{lerner-IMRN2012,MR3127380,MR4007575} and
Conde-Alonso and Rey~\cite{MR3521084} independently gave 
a new and simpler proof.  As part of their proofs
they introduced the technique of sparse domination, which use a
generalization of the classic dyadic cubes introduced by Calder\'on
and Zygmund.   A collection of cubes
\[ \D = \bigcup_{k\in \Z} \D_k \]
is a dyadic grid if for each $k\in \Z$, the cubes in $\D_k$ have
disjoint interiors, sidelengths $2^{-k}$, and their union is all of
$\rn$.  Further, given any two cubes $P,\,Q\in \D$, we must have that
$P\subset Q$, $Q\subset P$, or $P\cap Q = \emptyset$.   Examples of
such families include the translations of the standard dyadic grid
\[ \D^t = \{ 2^{-k}[0,1)^n + m +(-1)^kt : k \in \Z, m \in \Z^n \},
  \quad t \in \{0,\pm \tfrac{1}{3}\}^n.  \]
These collection share the property that given any cube $Q$ there exists
$t$ and $P \in \D^t$ such that $Q\subset P$ and $\ell(P) \leq
8\ell(Q)$.  (See~\cite{MR3092729,CruzUribe:2016ji}.)

Let
$\D$ be any dyadic grid.  A collection $\Ss\subset
\D$ is said to be sparse if there exist pairwise
disjoint sets $\{ E(Q) : Q \in \Ss \}$ and a constant $0<\eta<1$, such that for each $Q$, $E(Q)
\subset Q$ and $\eta |Q|\leq |E(Q)|$.   A sparse operator is an averaging
operator of the form
\[ T_{\Ss}f(x) = \sum_{Q\in \Ss} \left( \avgint_Q
  f(y)\,dy\right)\cdot \chi_Q(x). \]
Lerner and Nazarov and Conde-Alonso and Rey  showed that given a
singular integral $T$, for every bounded function of compact support there
exist a finite collection of dyadic grids $\{\D_k\}_{k=1}^N$ and
sparse families $\{\Ss_k\}_{k=1}^N$ such that
\[  |Tf(x)| \leq C \sum_{k=1}^N T_{\Ss_k}(|f|)(x).  \]
Given this, the $A_2$ conjecture reduces to proving the corresponding
estimates for sparse operators, which is much simpler.  This is done in two steps.  First,
weighted $L^2$ bounds are proved using an argument
that was essentially proved in~\cite{dcu-martell-perez}.  (See
also~\cite{dcu-paseky} where it is worked out in detail.)   Then $L^p(w)$ bounds follow from the
quantitative version of Rubio de Francia extrapolation.
Alternatively, a direct proof that holds for all $p>1$ was
given by Moen~\cite{MR3000426}.

\smallskip

There is another approach to weighted weak type inequalities that was
first investigated in the 1970s, but has recently gotten renewed
attention for its applications to matrix weights.  Given a weight
$w$, and an operator $S$, define  $S_w$ by
\[ S_wf(x) = w(x) S(w^{-1}f)(x).  \]
By a change of variables we have that
\begin{equation} \label{eqn:cov-operator}
  \|S_wf\|_{L^p} \leq C\|f\|_{L^p}
\end{equation}
if and only if
\[ \|Sf\|_{L^p(w)} \leq C\|f\|_{L^p(w)}. \]
By Chebyshev's inequality, this
implies the weak $(p,p)$ inequality
\begin{equation} \label{eqn:mult-weak1}
 \sup_{\lambda>0} \lambda \| \chi_{\{ x : |S_wf(x)| > \lambda \}} \|_{L^p(\rn)}
  \leq C \|f\|_{L^p(\rn)}.
\end{equation}
Inequalities of this form are referred to as multiplier weak
$(p,p)$ inequalities.   The following result is due to Muckenhoupt and
Wheeden~\cite{MR447956} on the real line, and
was extended to all $n\geq 1$ in~\cite{cruz-uribe-martell-perez05}.
Note that the only difficult case is $p=1$.

\begin{theorem} \label{thm:multiplier-weak}
  Given $1\leq p<\infty$ and $w\in \A_p$,
  inequality~\eqref{eqn:mult-weak1} holds when $S$ is either the
  Hardy-Littlewood maximal operator or a Calder\'on-Zygmund singular integral.
\end{theorem}

In~\cite{MR4859193} the authors proved a quantitative version of
Theorem~\ref{thm:multiplier-weak}.  Their results are only sharp when
$p=1$, as was proved by Lerner, {\em et al.}~\cite{MR4824922}.  Some
improvements of their results when $p>1$ were given by the same
authors in~\cite{lerner2024improvedweightedweaktype}.

\begin{remark}
  Surprisingly, Muckenhoupt and Wheeden showed that the $\A_p$
  condition is not necessary for multiplier weak type inequalities for
  either the maximal operator or the Hilbert transform, and that for
  $p>1$ they found different necessary conditions for these operators.  Recently,
  Sweeting~\cite{sweeting2024} proved that their necessary condition
  for the maximal operator is sufficient when $p>1$.
\end{remark}

\smallskip

There is a theory of off-diagonal scalar
weighted norm inequalities.  These are associated to the fractional
integral operator
\[  I_\alpha f(x) = \int_{\rn} \frac{f(y)}{|x-y|^{n-\alpha}}\,dy, \]
where $0<\alpha<n$, and the fractional maximal operator
\[ M_\alpha f(x) = \sup_Q |Q|^{\frac{\alpha}{n}} \avgint |f(y)|\,dy
  \cdot \chi_Q(x). \]
Muckenhoupt and Wheeden~\cite{muckenhoupt-wheeden74} introduced the
weight classes $\A_{p,q}$, where $1\leq p<\frac{n}{\alpha}$ and
$\frac{1}{p}-\frac{1}{q}=\frac{\alpha}{n}$, that generalize the
Muckenhoupt $\A_p$ classes:
\begin{equation} \label{eqn:Apq}
 [w]_{\A_{p,q}} = \sup_Q < \infty.  
\end{equation}
They proved the following result.

\begin{theorem} \label{thm:frac-op}
  Fix $0<\alpha<n$ and $1\leq p<\frac{n}{\alpha}$, and define $q$ by
  $\frac{1}{p}-\frac{1}{q}=\frac{\alpha}{n}$.  Given a weight $w$, if $S_\alpha$ is either
  the fractional integral operator $I_\alpha$ or the fractional
  maximal operator $M_\alpha$, the following are equivalent:
  \begin{enumerate}

  \item $w\in \A_{p,q}$;
  \item $S_w$ satisfies the weak $(p,q)$ inequality
    \[  \sup_{\lambda>0} \lambda \| \chi_{ \{x\in \rn : |S_\alpha
        f(x)|>\lambda\} } \|_{L^q(w)}
        \leq 
        C\|f\|_{L^p(w)};\]
      \item additionally,  if $p>1$, $S_\alpha$ satisfies the strong $(p,q)$
        inequality
        \[ \|S_\alpha f\|_{L^q(w)} \leq C\|f\|_{L^p(w)}. \]
      \end{enumerate}
    \end{theorem}
    
The original proof of Theorem~\ref{thm:frac-op} by  Muckenhoupt and Wheeden was to prove norm
inequalities for $M_\alpha$, and then to use a good-$\lambda$
inequality similar to~\eqref{eqn:good-lambda-sio}  to prove
norm inequalities for $I_\alpha$.  This approach only gives
qualitative estimates for the constant.  The sharp constant in terms
of the $[w]_{\A_{p,q}}$ characteristic was proved by Lacey, {\em et
  al.}~\cite{MR2652182}.   A refined version of their proofs that
developed the theory of off-diagonal dyadic sparse operators, was given
in~\cite{DCU-KMfrac2,MR3224572}; see also~\cite{CruzUribe:2016ji}.
Much of this was based on the seminal work of P\'erez~\cite{perez94}.

\section{Jones factorization and Rubio de Francia
  extrapolation}
\label{section:jones-rubio}

As we noted above, the proof that sparse operators satisfy weighted
$L^p$ bounds with the conjectured sharp  constant follows at once from
the case $p=2$ by the Rubio
de Francia extrapolation theorem.  This theorem has a long history and
we refer the reader to~\cite{MR2797562} for complete details.
Here we discuss a version of theorem that has several
features not present in the original.  First, we give a
 version that gives quantitative estimates for the constants.  This
 was first proved when $p_0<\infty$ by~ Dragi\v{c}evi\'{c}, {\em et
   al.}~\cite{MR2140200}.   A simpler proof is given in
 Duoandikoetxea~\cite{10.1016/j.jfa.2010.12.015} and also
 in~\cite{MR2797562}.

 \begin{remark}
   This is often referred to as the ``sharp constant'' version of
   extrapolation because when applied to singular integrals it gives
   the sharp constants for all $p$ when starting from $p=2$, and so
   allowed proof of the $A_2$ conjecture.  Similarly, it gives the
   sharp constants for the dyadic square function:
   see~\cite{dcu-martell-perez}.  However, as we will see below for
   matrix weights, it may not  give the best constant, even when starting
   from an inequality with the sharp constant.
 \end{remark}
 
 Second, taking advantage of the uniform definition of $\A_p$, $1\leq
 p \leq \infty$, we give
 a version that lets us start from inequalities when $p=\infty$.  A qualitative version of
this result was proved by Harboure, {\em et al.}~\cite{MR944321}; they
implicitly used the class of $\A_\infty$ weights.  The quantitative version
of this result was proved by Nieraeth~\cite{MR4000248} and
independently in~\cite{mb-2022}.  

Third, we want to state an abstract version of extrapolation.  Rather
than working with a specific operator (which actually plays no direct role in
the proof), we work with pairs of nonnegative functions $(f,g) \in
\F$, which we will refer to as extrapolation pairs.  Hereafter, we
will implicitly assume that given an inequality of the form 
\begin{equation} \label{eqn:extrapol-pairs}
\|f\|_{X} \leq C\|g\|_{X}, \qquad (f,g)\in \F,  
\end{equation}
where $X$ is some Banach function space, we mean that this inequality holds for all pairs $(f,g)\in \F$ such
that $0< \|f\|_X, \, \|g\|_X < \infty$.  
This approach to extrapolation was introduced
in~\cite{cruz-uribe-perez00} and systematically developed
in~\cite{cruz-uribe-martell-perez04, MR2797562}.  Intuitively, to
prove bounds for an operator $T$, one takes the pairs $(|Tf|,|f|)$;
for a discussion of the technical details of this approach,
see~\cite{MR2797562,dcu-paseky}. 

\begin{theorem} \label{thm:rubio}
Given $p_0$, $1\leq p_0\leq \infty$, suppose that for a family of
extrapolation pairs $\F$
and 
for all $w_0\in \A_{p_0}$, the inequality
\[ \|f\|_{L^{p_0}(w_0)} 
  \leq
  N_{p_0}([w_0]_{\A_{p_0}}) \|g\|_{L^{p_0}(w_0)}, \qquad (f,g)\in \F, \]
holds.  Then for all $p$, $1<p<\infty$,  and all $w\in \A_p$, 
\[ \|f\|_{L^{p}(w)} 
  \leq
  N_{p}([w]_{\A_{p}}) \|g\|_{L^{p}(w)}, \qquad (f,g)\in \F, \]
where 
\[  N_p ([w]_{\A_{p}}) 
= 
C(n,p_0,p)N_{p_0}\big(c(n,p_0,p) [w]_{\A_{p}} ^{\max\{\frac{p}{p_0}, \frac{p'}{p_0'}\}}\big).\]
\end{theorem}

\begin{remark}
  One feature of this result is that it does not allow extrapolation
  to the endpoints:  starting with $1<p_0<\infty$, it is not possible
  to get to either $p=1$ or $p=\infty$.  
\end{remark}

\smallskip

We will not prove Theorem~\ref{thm:rubio} here.  For complete proofs
(when $p<\infty$) see~\cite{MR2797562}; for a discussion of the
essential ideas of the proof, see~\cite{dcu-paseky}.   Here we want to
concentrate on the key ideas underlying the proof.  

The proof  ultimately depends on four things.  First, it requires the
``duality'' implicit in the definition of $\A_p$ weights:  that for
$1\leq p \leq \infty$, $w\in \A_p$ if and only if $w^{-1}\in \A_{p'}$.  
Second, it requires the quantitative bound for the strong $(p,p)$
inequality for the maximal operator discussed above: for
    $1<p \leq \infty$,
    $ \|Mf\|_{L^p(w)} \leq C[w]_{A_p}^{p'}\|f\|_{L^p(w)}$.
    Third, it requires the Jones factorization theorem.

    \begin{theorem}  \label{thm:jones}
Given $1<p<\infty$ and a weight $w$, $w\in \A_p$
   if and only if there exist $w_0\in \A_1$ and $w_1 \in \A_\infty$
   such that
   \[ w=w_0^{\frac{1}{p}}w_1^{\frac{1}{p'}}.  \]
 \end{theorem}

 Theorem~\ref{thm:jones} was first proved by Jones~\cite{Jones},
   but a much more elementary proof was given by Coifman, Jones, and
   Rubio de Francia~\cite{MR687639}.  (See also~\cite{dcu-paseky}.)
   In fact, in the proof of extrapolation, we only use the easier
  half of the Jones factorization theorem:  that if $w_0\in A_1$ and
  $w_1 \in \A_\infty$, then
  $w=w_0^{\frac{1}{p}}w_1^{\frac{1}{p'}}\in \A_p$. This property,
  which we refer to as
  ``reverse factorization,'' 
  follows at once from the definition of the $\A_p$ classes.  

Fourth, we need a means of constructing $\A_p$ weights, which, by
applying Theorem~\ref{thm:jones}, means we need to construct $\A_1$ and
$\A_\infty$ weights (which are just inverses of $\A_1$ weights).  To
do so we use the Rubio de Francia iteration
operator, which is also fundamental to the proof of the harder
implication in the Jones
factorization theorem.  
Fix $w\in \A_p$,
$1<p<\infty$.  Given a
nonnegative function $h$,  define
\[ \Rdf h(x) = \sum_{k=0}^\infty \frac{M^kh(x)}{2^k
    \|M\|_{L^p(w)}^k}, \]
where $\|M\|_{L^p(w)}$ is the operator norm of $M$ on $L^p(w)$.
It then follows from the definition and the properties of $\A_p$ weights that
\begin{enumerate}

\item $h(x) \leq \Rdf h(x)$;

\item $\|\Rdf h\|_{L^p(w)} \leq 2 \| h\|_{L^p(w)}$;

  \item $\Rdf h \in A_1$ and $[\Rdf h]_{A_1} \leq 2\|M\|_{L^p(w)} \leq
    C[w]_{A_p}^{p'}$.

  \end{enumerate}
  Each of these estimates plays an important role in the proof.
  
  \smallskip

  We briefly consider three extensions of extrapolation.
  First, one consequence of working with abstract extrapolation pairs
  is that it immediately lets us prove vector-valued inequalities of
  the form
  \begin{equation} \label{eqn:scalar-vv}
 \bigg\|  \bigg(\sum_{k=1}^\infty |wf_k|^q
    \bigg)^{\frac{1}{q}}  \bigg\|_{L^p(\rn)}
    \leq
    C  \bigg\|  \bigg(\sum_{k=1}^\infty |wg_k|^q
    \bigg)^{\frac{1}{q}}  \bigg\|_{L^p(\rn)}
  \end{equation}
  when $1< p,\,q<\infty$, $w\in \A_p$, and $\{(f_k,g_k)\}_{k=1}^\infty
  \subset \F$.  It also lets us apply
  extrapolation to prove (classical) weak type inequalities
\[ \sup_{\lambda>0} \lambda \| \chi_{\{ x : |f(x)| > \lambda \}}\|_{L^p(w)}
  \leq
  C\|g\|_{L^p(w)}, \qquad (f,g) \in \F. \]
For details, see~\cite{MR2797562}.

\begin{remark}
It is an open question whether there exists a version of Rubio de
Francia extrapolation which can be applied to the multiplier weak type
inequalities.
\end{remark}

Second, we can apply extrapolation to more general classes of
weights.  By a basis, we mean a nonempty collection of open, bounded sets $\B$
in $\rn$.  We can then define $\A_p$ classes with respect to a basis:
for $1\leq p\leq \infty$, we say $w\in \A_{p,\B}$ if
\begin{equation} \label{eqn:basis-Ap}
 [w]_{\A_{p,\B}} = \sup_{B\in \B} |B|^{-1}
  \|w\chi_B\|_{L^p}\|w^{-1}\chi_B\|_{L^{p'}} < \infty.  
\end{equation}
We say that a basis $\B$ is a Muckenhoupt basis if the maximal
operator
\begin{equation} \label{eqn:basis-max}
 M_{\B}f(x) = \sup_{B\in \B} \avgint_B |f(y)|\,dy \cdot \chi_B(x) 
\end{equation}
satisfies $\|M_{\B}f\|_{L^p(w)} \leq C\|f\|_{L^p(w)}$, $1<p\leq
\infty$, whenever $w\in \A_{p,\B}$. There are several important examples of Muckenhoupt bases
besides the basis of cubes $\Q$:
\begin{itemize}
    \item the set $\D$ of dyadic cubes (this is implicit in the
      proofs of bounds for the Hardy-Littlewood maximal operator);
    \item the set $\Rdf$ of rectangles with sides parallel to the coordinate axes (e.g., the strong maximal operator);
    \item more generally, the multiparameter basis $\Rdf^\alpha$,
      $\alpha=(a_1,\ldots,a_j)$, consisting of sets of the form
      $Q_1\times\cdots \times Q_j$, where $Q_i$ is a cube in
      $\R^{a_i}$ (see Garc\'\i a-Cuerva and Rubio de Francia~\cite{garcia-cuerva-rubiodefrancia85});
    \item the Zygmund basis $\Zz$ of rectangles in $\R^3$ with sides
      parallel to the coordinate axes and whose sidelengths satisfy
      $(s,t,st)$ for some $s,\,t>0$ (see Fefferman and Pipher~\cite{MR1439553}). 
\end{itemize} 
With this definition, the proofs of extrapolation (e.g., of Theorem~\ref{thm:rubio}) go through without
any essential changes.  This is done in great detail
in~\cite{MR2797562}; note that the above fundamental ideas from the
proof carry through immediately for Muckenhoupt bases.

\smallskip

Finally, we note that extrapolation can be extended to the
off-diagonal setting, using the $\A_{p,q}$ weights \eqref{eqn:Apq}  of Muckenhoupt and
Wheeden,  A qualitative version of the following result was proved
(when $q_0<\infty$) by Harboure, {\em et al.}~\cite{MR944321}.
A quantitative version (again when $q_0<\infty$), which yielded the sharp constant estimates for
fractional integrals, was proved by Lacey, {\em et
  al.}~\cite{MR2652182}.  Their proof generalized an older proof of
extrapolation that depended on a technical lemma due go Garc\'\i
a-Cuerva~\cite{MR684631}. (See
also~\cite{garcia-cuerva-rubiodefrancia85}.)  A different, simpler
proof, which does not get the best possible constant, is in~\cite{MR2797562}.  (For the
explicit computation of this constant, in the setting of matrix
weights, see~\cite{dcu-fs-AFM2025}.)   For a simpler proof of the
``sharp constant'' result of Lacy, {\em et al.}, see
Duoandikoetxea~\cite{10.1016/j.jfa.2010.12.015}.  For another proof, including the case
$q_0=\infty$, see the forthcoming paper~\cite{dcu-fs-2025}.

\begin{theorem} \label{thm:rubio-offdiag}
Fix $s$, $1\leq s<\infty$. Given $p_0,\, q_0$, $1\leq p_0\leq q_0 \leq
\infty$, $\frac{1}{p_0}-\frac{1}{q_0}= \frac{1}{s'}$, suppose that for
some family of extrapolation pairs $\F$
and 
for all $w_0\in \A_{p_0,q_0}$, the inequality
\[ \|f\|_{L^{q_0}(w_0)} 
  \leq
  N_{p_0,q_0}([w_0]_{\A_{p_0,q_0}}) \|g\|_{L^{p_0}(w_0)}, \qquad (f,g)\in \F, \]
holds.  Then for all $p,\,q$, $1<p\leq q<\infty$,
$\frac{1}{p}-\frac{1}{q}= \frac{1}{s'}$,  and all $w\in \A_{p,q}$, 
\[ \|f\|_{L^{q}(w)} 
  \leq
  N_{p,q}([w]_{\A_{p,q}}) \|g\|_{L^{p}(w)}, \qquad (f,g)\in \F, \]
where 
\[  N_{p,q} ([w]_{\A_{p,q}}) 
= 
C(n,p_0,p,s)N_{p_0}\big(c(n,p_0,p,s) [w]_{\A_{p,q}} ^{\max\{\frac{q}{q_0}, \frac{p'}{p_0'}\}}\big).\]
\end{theorem}

In the proof of Theorem~\ref{thm:rubio-offdiag}, a key fact is that
$w\in \A_{p,q}$ if and only if $w^s\in \A_r$.  This equivalence follows immediately
from the definitions.

\section{Matrix weights and matrix weighted norm inequalities}
\label{section:matrix}

We now turn to the problem of generalizing weighted norm inequalities
to the setting of vector-valued functions and matrix weights.  We
first  define some notation.  Let
$\vecf=(f_1,\ldots,f_d)$ be a measurable $\rd$-valued function.  Given a singular integral $T$, define  it
acting on $\vecf$ by
\[ T\vecf = (Tf_1,\ldots,Tf_d).  \]
It is immediate that if $1<p<\infty$ and $\vecf \in L^p(\rn,\rd )$,
then $\|T\vecf\|_{ L^p(\rn,\rd )} \leq C \|\vecf\|_{ L^p(\rn,\rd )}$.

To define matrix weights, recall that $\Sd$ is the set of $d\times d$,
self-adjoint, positive semidefinite matrices.  A matrix weight is a
measurable function  $W : \rn \rightarrow \Sd$.  We will assume that
the entries of $W$ are finite and that $W$ is invertible almost
everywhere.  (This assumption is in some sense necessary:
see~\cite{dcu-fs-AFM2025}.)  Define a scalar weight using the operator norm of $W$:
\[ |W(x)|_{\op}
  =
  \sup_{\xi \in \rd, |\xi|=1} |W(x)\xi|. \]
For $1\leq p \leq \infty$, we define the matrix weighted space
$L^p(W)=L^p(W,\rn,\rd )$ by
$\|\vecf\|_{L^p(W)} = \|W\vecf\|_{L^p(\rn,\rd)}$.  When $p<\infty$,
this becomes
\[ \|\vecf\|_{L^p(W)} = \bigg( \int_\rn
  |W(x)\vecf(x)|^p\,dx\bigg)^{\frac1p}. \]
Note that when $d=1$, this reduces to the scalar space $L^p(w)$.

\begin{remark}
  Originally, these spaces were defined for $p<\infty$ with the matrix $W$ replaced
  by $W^{\frac1p}$, so that when $d=1$ it would reduce to the
  classical definition of weighted norm inequalities.  Note that $W^{\frac1p}$ is well-defined since $W$ is positive
semidefinite.   For a variety of reasons we believe this definition is
a better one.  In particular, see the discussion below on left openness.
\end{remark}

With this notation, the problem first considered by Nazarov, Treil and
Volberg can be stated as follows:  given $1<p<\infty$, prove there is a Muckenhoupt type
condition on matrix weights so that the inequality
\[\|T\vecf\|_{L^p(W)} \leq \|\vecf\|_{L^p(W)}.   \]
holds for singular integrals.    Treil and Volberg~\cite{MR1428818}  first
solved this problem on the real line for the Hilbert transform  when
$p=2$.  They showed that this inequality holds if $W$ satisfies the
following  analog of the $A_2$
condition:
\begin{equation} \label{eqn:matrix-A2}
  \sup_Q \bigg| \left(\avgint_Q W^2(x)\,dx\right)^{\frac{1}{2}}
  \left(\avgint_Q W^{-2}(x)\,dx\right)^{\frac{1}{2}}\bigg|_{\op}
  <\infty. 
\end{equation}
This condition, however, cannot be extended to the case $p\neq 2$.

Ultimately, three equivalent definitions were developed for $\A_p$,
$1\leq p\leq \infty$, each with its
own strengths and weaknesses.  The first definition of matrix $\A_p$
was in terms of
norm functions.  It  was
conjectured by Treil~\cite{MR1030053} and used by Nazarov and Treil~\cite{MR1428988}
and by Volberg~\cite{MR1423034} to prove matrix weighted norm
inequalities for the Hilbert transform.  Their idea was to replace
matrices with norms on $\rd $.  Here we sketch their definition; for
complete details see the above references or~\cite{mb-2022}.

Recall that a function $\rho : \R^d \rightarrow [0,\infty)$ is a norm
if,
given
$v,\,w\in \rd$ and $\alpha \in \R$:  
$\rho(v)=0$ if and only if $v=0$;
$\rho(v+w) \leq \rho(v)+\rho(w)$;
and $\rho(\alpha v)=|\alpha|\rho(v)$. 
Define the dual of a norm $\rho$ to be the norm $\rho^*$ given by 
\[ \rho^*(v) = \sup_{w\in \rd , \rho(w)\leq 1} |\langle v,
  w\rangle|. \]
Let  $\rho : \rn \times \rd 
\rightarrow [0,\infty)$ be a measurable function such that for
a.e. $x\in \rn$, $\rho(x,\cdot)$ is a norm on $\rd $.
For instance, given a matrix weight $W$, we can define a norm function
$\rho_W$ by $\rho_W(x,v)=|W(x)v|$. 
Finally, given a cube $Q$, and $1\leq p<\infty$, define the average of a
norm function on $Q$ by
\[ \langle \rho \rangle_{p,Q}(v)
  =
  \| \rho(\cdot,v)\|_{p,Q}
  =
  \bigg(\avgint_Q \rho(x,v)^p\,dx\bigg)^{\frac1p}. \]
When $p=\infty$, let $ \langle \rho \rangle_{\infty,Q}(v)=
\|\rho(\cdot,v)\chi_Q\|_{L^\infty(\rn)}$.  
 For all
$p$, $1\leq p \leq \infty$, we have that  for
every $v\in \rd $,
\[  \langle \rho \rangle_{p,Q}^* (v) \leq  \langle \rho^*
  \rangle_{p',Q}(v).  \]
We define  $\rho$ to be in $\A_p$ if the reverse inequality holds up
to a constant independent of $v$:
\[  \langle \rho^* \rangle_{p',Q}(v)
  \leq
  C  \langle \rho \rangle_{p,Q}^* (v) \]
When $d=1$, if $\rho=\rho_w$ for some scalar weight $w$, then this
reduces to the scalar $\A_p$ condition after rearranging terms.

Given any norm function we can associate to it a norm function induced
by a matrix that is equivalent to it. Let $\rho$ be a norm, and define
$K_\rho$ to be the closed unit ball in $\R^d$ with respect to this norm:
\begin{equation} \label{eqn:norm-ball}
K_\rho = \{ v \in \rd : \rho(v) \leq 1 \}. 
\end{equation}
Then $K_\rho$ is a convex set in $\rd$; associated to it is its John
ellipsoid:  a unique ellipsoid $E$ of maximum volume such that
$E \subset K_\rho \subset \sqrt{d}E$.  We can write $E=W\mathbf{B}$,
where $W \in \Ss_d$.  As a consequence, $\rho(v) \approx \rho_W(v)$
for all $v\in \rd$.  Moreover, if $\rho$ is a norm function, then we
can choose the associated matrix function $W$ to be measurable.
(See~\cite[Theorem~4.11]{mb-2022}.)   Consequently, it suffices to
restrict our attention to norms induced by matrices.  

Working with matrices,
Roudenko~\cite{MR1928089} gave an equivalent definition of 
$\A_p$ that more closely resembled the scalar definition:  for
$1<p<\infty$, $W\in \A_p$ if and only if
\begin{equation} \label{eqn:roudenko1}
 [W]_{\A_p} = \sup_Q \avgint_Q \bigg( \avgint_Q
    |W^{\frac{1}{p}}(x)W^{-\frac{1}{p}}(y)|_{\op}^{p'}\,dy\bigg)^{\frac{p}{p'}}\,dx <
    \infty. 
  \end{equation}
  Frazier and Roudenko~\cite{MR2104276} also introduced the concept of matrix $\A_1$
  weights:  $W\in \A_1$ if
  \begin{equation} \label{eqn:roudenko2}
  [W]_{\A_1} = \sup_Q \esssup_{x\in Q}
    \avgint_Q |W^{-1}(x)W(y)|_{\op} \,dy < \infty. 
  \end{equation}
  The definition of $\A_\infty$ was given in~\cite{mb-2022}:  $W\in \A_\infty$ if
  \begin{equation} \label{eqn:roudenko3}
  [W]_{\A_\infty} = \sup_Q \esssup_{x\in Q}
    \avgint_Q |W(x)W^{-1}(y)|_{\op} \,dy < \infty. 
  \end{equation}
  For the equivalence of this definition with the previous one, see~\cite{mb-2022}.
  One advantage of this definition is that it is straightforward to
  show that $W\in \A_p$, $1\leq p\leq \infty$, if and only if the
  averaging operator $A_Q\vecf(x) = \avgint_Q \vecf(y)\,dy\cdot
  \chi_Q(x)$, is bounded on $L^p(W)$.  This in turn can be used to
  prove norm inequalities for convolution operators and norm
  convergence of approximate identities.
  (See~\cite{MR3544941}.)  When $d=1$, this definition immediately
reduces to the definition for scalar weights given above.  

  The third equivalent definition of matrix $\A_p$ weights is in terms
  of reducing operators.  This approach was introduced by
  Volberg~\cite{MR1423034} and further developed and used by
  Goldberg~\cite{MR2015733}. Given a matrix $W \in \Ss_d$ and a cube
  $Q$, for $1\leq p \leq \infty$ we define the norm $\langle \rho_W
  \rangle_{p,Q}$.  As noted above, there exists a matrix $\W^p_Q$
  which induces an equivalent norm.  We call this matrix the reducing
 operator associated to $W$ on $Q$.  Similarly, we let $\overline{W}^{p'}_Q$
  be the reducing operator associated with $\langle \rho_{W^{-1}}
    \rangle_{p',Q}$.   Then $W \in \A_p$ if and only if
    \[ [W]_{\A_p}^R = \sup_Q |\W^p_Q \overline{W}^{p'}_Q|_{\op} < \infty. \]
  Moreover, we have that $ [W]_{\A_p}^R \approx  [W]_{\A_p}$.   (Again,  see~\cite{mb-2022}.) 
The
reducing operator definition is very useful in technical estimates
related to the maximal operator:  see below.  Here we note that, using
the fact that if $A$ and $B$ are self adjoint matrices, then
$|AB|_{\op}=|BA|_{\op}$, we have that $W\in \A_p$ if and only if
$W^{-1} \in \A_{p'}$ and $[W^{-1}]_{\A_{p'}}\approx [W]_{\A_p}$.  Using \eqref{eqn:roudenko1} above, this
equivalence only follows when $p=2$.  

\smallskip

Norm inequalities for   Calder\'on-Zygmund singular
integrals in $\R^n$ were proved by Christ and
Goldberg~\cite{MR1813604,MR2015733} for $1<p<\infty$.  Reducing
operators were central to their proof.    Treil and Volberg had earlier noted that a key
obstacle to proving matrix weighted inequalities was the lack of a
maximal operator that did not lose the geometric
information embedded in a vector-valued function.  A key
component of the proofs in~\cite{MR1813604,MR2015733} is a scalar-valued, matrix
weighted maximal operator, now referred to as the Christ-Goldberg
maximal operator:
\[ M_W{\vecf}(x) =
  \sup_Q \avgint_Q |W(x)W^{-1}(y){\vecf}(y)|\,dy\cdot \chi_Q(x).  \]
The motivation for this definition comes from change of variables
introduced in the scalar case when defining multiplier weak type
inequalities.  Note that when $d=1$, $M_wf(x)=w(x)M(w^{-1}f)(x)$, so
the boundedness of $M_w$ on $L^p(\R^n)$ is equivalent to the
boundedness of the Hardy-Littlewood maximal operator on $L^p(w)$.   
We can make a similar definition for singular integral operators;
since they are linear the matrix $W(x)$ can be pulled in or out of the operator:
\[ T_W\vecf(x) = W(x)^{\frac1p}T(W^{-\frac1p}\vecf)(x). \]
The following theorem then holds.

\begin{theorem} \label{thm:goldberg}
  Given $1<p<\infty$ and a matrix weight $W \in \A_p$, then
  \begin{enumerate}

  \item The maximal operator $M_W$ is bounded from $L^p(\rn,\rd)$ to $L^p(\rn)$;

    \item If $T$ is a Calder\'on-Zygmund singular integral, then $T_W$
      is bounded on $L^p(\rn,\rd)$.

    \end{enumerate}
    Moreover, if $M_W$ is bounded, or if $T$ is bounded
    and 
    nondegenerate, then $W\in \A_p$. 
  \end{theorem}

  \begin{remark}
    When $p=\infty$, $M_W$ is bounded on $L^\infty(\rn)$ if and only
    if  $W\in \A_\infty$.  See~\cite{mb-2022}.
  \end{remark}

The proof of Theorem~\ref{thm:goldberg}  roughly follows the
proofs of the corresponding results in the scalar case, but a number
of technical obstacles arise.
  The first step was to prove that if $W\in \A_p$, $1<p<\infty$, then $M_W$ is bounded
from  $L^p(\rn,\rd)$ into $L^p(\rn)$.  This required the introduction
of an auxiliary maximal operator $M_W'$, defined by
\[ M_W'{\vecf}(x) =
  \sup_Q \avgint_Q |\W_Q^p W^{-1}(y){\vecf}(y)|\,dy\cdot \chi_Q(x).  \]
The proof that the auxiliary maximal operator maps $L^p(\rn,\rd)$ into
$L^p(\rn)$ is done by first  using a covering lemma argument (e.g.,
the Vitali covering lemma or Calder\'on-Zygmund cubes) to show
that it satisfies a weak $(p,p)$ inequality.  In fact, by an argument
which is analogous to using the reverse H\"older inequality in the
scalar case to prove
that if $w\in \A_p$, then $w\in \A_{p\pm \epsilon}$,
they showed that there
exists $\epsilon>0$ such that it satisfies weak $(p\pm \epsilon, p\pm
\epsilon)$ inequalities.  The desired strong type inequality then
follows by Marcinkiewicz interpolation.
Finally,  while there is no direct pointwise relationship between
$M_W$ and $M_W'$, it follows from  a  stopping time
argument that
\[ \|M_W\vecf\|_{L^p(\rn,\rd)} \leq C \|M_W'\vecf\|_{L^p(\rn,\rd)}. \]

Given these bounds for the maximal operators $M_W$ and $M_W'$, the
fact  that singular integrals are bounded is gotten by adapting the ideas of Coifman and Fefferman to
prove a good-$\lambda$ inequality.  Define
\[ T_W^*\vecf(x) = \sup_{\epsilon>0}
  |W^{\frac1p}(x)T_\epsilon(W^{-\frac1p}\vecf)(x)|. \]
Then for every smooth function $\vecf$ with compact support, they
proved that there exist constants $0<b<1$ and $c>0$ such that for all
$\lambda>0$,
\begin{multline*}
  |\{x\in \rn : T_W^*\vecf(x)>\lambda,
  \max\{ M_W'\vecf(x), M_W\vecf(x)\}<c\lambda\}| \\
  \leq
  \frac{1}{2}b^p|\{x\in \rn : T_W^*\vecf(x) > b\lambda\}|. 
\end{multline*}
Bounds for $T^*$ and so for $T$ then follow by a standard
level-set argument.  

\smallskip

As in the scalar case, this approach to proving matrix weighted norm
inequalities does not yield quantitative estimates on the constant in
terms of the $[W]_{\A_p}$ characteristic.  After the proof of the
  scalar $A_2$ conjecture by Hyt\"onen,  it was conjectured that the
  corresponding result holds for matrix weights:  for $1<p<\infty$, 
  \[ \|T\vecf\|_{L^p(W)}
    \leq
    C[W]_{\A_p}^{\max\{p,p'\}}\|\vecf\|_{L^p(W)}.  \]
  This problem was considered by Bickel, Petermichl and
Wick~\cite{MR3452715}, Pott and Stoica~\cite{MR3698161} and Culiuc, Di
Plinio and Ou~\cite{MR3818613} when
$p=2$.    Here we will consider the approach of  Nazarov, Petermichl, Treil and
Volberg~\cite{MR3689742}.  They showed that when $p=2$,
  \[ \|T\vecf\|_{L^2(W)}
    \leq
    C[W]_{\A_2}^{3}\|\vecf\|_{L^2(W)}.  \]
Their proof is based on a deep
generalization of the sparse domination estimates described above in
the scalar case. 

The sparse operator introduced in~\cite{MR3689742} replaces
vector-valued functions by averages that are convex sets.  They show
that given $\vecf \in L_c^\infty(\rn,\rd)$, there exists a finite
collection of dyadic grids $\{\D_k\}_{k=1}^N$ and sparse sets  $\{\Ss_k\}_{k=1}^N$ such that
\begin{equation} \label{eqn:nptv1}
  T{\vecf}(x) \in
C\sum_{k=1}^N T_{\Ss_k}\vecf(x) 
  = C\sum_{k=1}^N\sum_{Q\in \Ss_k}
  \langle\langle{\vecf}\rangle\rangle_Q \chi_Q(x), 
\end{equation}
where $\langle\langle{f}\rangle\rangle_Q$ is the convex set
\[  \langle\langle{\vecf}\rangle\rangle_Q
  = \bigg\{ \avgint_Q k(y) {\vecf}(y)\,dy : k \in L^\infty(Q),
  \|k\|_{L^\infty(Q)} \leq 1 \bigg\},  \]
and the sum is the (infinite) Minkowski
sum of convex sets.   They referred to this estimate as a convex body
sparse domination.  To complete their proof, however, they did not
work directly with these convex set-valued functions, but rather
replaced them by vector-valued sparse operators of the form
\[ \tilde{T}_\Ss {\vecf}(x) = \sum_{Q\in \Ss} \avgint_Q
  \varphi_Q(x,y){\vecf}(y)\,dy, \]
where for each $Q$, $\varphi_Q$ is a real-valued function supported on
$Q\times Q$  such that, for each $x$,
$\|\varphi_Q(x,\cdot)\|_\infty \leq 1$.
These they estimated using square functions.

When $1<p<\infty$,  quantitative
results were proved in~\cite{MR3803292}.  They showed that
  \[ \|T\vecf\|_{L^p(W)}
    \leq
    C[W]_{\A_p}^{p+p'-1}\|\vecf\|_{L^p(W)}.  \]
They used the sparse domination result of Nazarov, {\em et al.} and
reduced to the vector-valued sparse operators, but
instead of using square function estimates, they adapted techniques
from the study of scalar, two weight
norm inequalities, introducing matrix analogs of the so-called $\A_p$
bump conditions.  Consequently, they also proved two weight
inequalities for singular integrals, the first
such results in the matrix weight setting.

This proof, while quantitative, did not follow the approach which
yielded sharp inequalities in the scalar case for sparse operators.
They used Orlicz maximal function estimates which, in the scalar
one weight case, were known to not produce the sharp estimate.
Initially, there was
hope that this proof could be refined to prove the $\A_2$ conjecture
for matrix weights.  In the scalar case the key tool is the 
``universal'' dyadic maximal operator.  Given any weight $w$ and a
dyadic grid $\D$, define the
dyadic operator $M_w^\D$ by
\begin{equation} \label{eqn:universal-max}
 M_w^\D f(x) = \sup_{Q\in \D} \frac{1}{w(Q)}\int_Q |f(y)|\,w(y)dy \cdot
  \chi_Q(x).  
\end{equation}
Dyadic covering arguments yield the estimate
\[ \bigg(\int_\rn M_w^\D f(x)^p\, w(x)dx \bigg)^{\frac1p}
  \leq p' \bigg(\int_\rn |f(x)|^p\, w(x)dx \bigg)^{\frac1p}. \]
(See~\cite{MR3000426}. While it is more natural to write this
inequality with the weight acting as a measure, by replacing $w$ by
$w^p$ we can write it with the weight acting as a multiplier.)   The
goal was to find an analogous operator
in the matrix case and adapt the scalar proofs to use it.
There is a natural candidate, based on the convex set-valued maximal
operator (see below), but it was shown to not be bounded on
$L^p(\rn,\rd)$:  see\cite{MR4491249}.  Moreover, very recently the exponent $3$
was shown to be the best possible when $p=2$ by Domelevo, Petermichl, Treil, and
Volberg~\cite{DPTV-2024}.  This bound is also the best possible for
sparse operators: see Treil and Volberg~\cite{tv-2025}.

\begin{remark}
  It is
an open question whether the exponent $p+p'-1$ is the best possible
for all $p$. We are tempted to conjecture that it is, but there is no
evidence either way.
\end{remark}

\smallskip

For many years, only strong type inequalities were considered in the
matrix case, as it did not seem possible to generalize the classical
weak type inequalities to use matrix weights.  But then,
in~\cite{MR4269407} the authors utilized the model of multiplier weak
type inequalities for scalar weights to define and prove the analogous
estimates for matrix weights.  The following result was proved in~\cite{MR4269407}  when
$p=1$, and for $1<p<\infty$ in~\cite{MR4859193}.

\begin{theorem} \label{thm:sio-weak-matrix}
  Fix $1\leq p<\infty$, and $W\in \A_p$, if $T$ is a
  Calder\'on-Zygmund singular integral operator, then
  \begin{equation} \label{eqn:swm1}
 \sup_{\lambda>0} \lambda |\{ x\in \rn :
    |T_W\vecf(x)|>\lambda\}|^{\frac1p}
    \leq C[W]_{\A_p}^{p+1} \|\vecf \|_{L^p(\rn)}.  
  \end{equation}
  Moreover, the same inequality holds if $T_W$ is replaced by the
  Christ-Goldberg maximal operator $M_W$. 
\end{theorem}

\begin{remark}
  As we noted above, the weak $(p,p)$ inequality follows from the
  strong $(p,p)$ inequality by Chebyshev's inequality.  The constant
  above is better in the range $1<p< 2$ when compared to the values
  gotten from the best known constant in the strong type inequality.
This exponent is sharp when $p=1$ (see \cite{MR4824922}), but it is not known what the best
constant is otherwise.  For some partial results,
see~\cite{lerner2024improvedweightedweaktype}. 
\end{remark}

\smallskip

The proof of Theorem~\ref{thm:sio-weak-matrix} uses convex body
sparse domination to reduce the problem to an estimate for a scalar
weighted sparse operator.  The proof then follows the same pattern as
the proof of the scalar version of this result.  We note that this
part of the proof  uses a clever restatement of the weak type
norm on the lefthand side:

    \[
        \|w T_\mathcal{S}(w^{-1}\,\cdot)\|_{L^p \rightarrow
          L^{p,\infty}}
        \approx \sup_{\|f\|_{L^p} = 1}\,\sup_{0<|E|<\infty}\,
        \inf_{\substack{E'\subseteq E \\ |E| \leq 2
            |E'|}}|E|^{\frac{1}{p}-1}|\langle w T_\mathcal{S}(w^{-1}f),\chi_{E'}\rangle|.
    \]
(See~\cite[Exercise~1.4.14]{grafakos08a}; this approach to estimating
weak type norms was previously used by Frey and Nieraeth~\cite{MR3897012}.)

\smallskip

As an application of
Theorem~\ref{thm:sio-weak-matrix} we describe another proof of the boundedness of
singular integrals in the matrix weighted setting that makes use of a
recently discovered structural property of matrix $\A_p$ weights.   As we noted
above, Bownik~\cite{Bow} showed that for the classical definition
of matrix $A_p$ weights, left openness does not hold.  More precisely,
he constructed
a matrix weight $W\in A_2$ that is not in $A_p$ for any $p<2$.
However, this property does hold using the multiplier definition given
above.  In~\cite{dcu-mp-2025} it was shown that if $W\in \A_p$,
$1<p<\infty$, then there exist $0<r<1<s$ such that $W\in \A_{rp}$ and
$W\in \A_{sp}$.   Using this fact, Penrod and
Sweeting~\cite{MP-BS-2026} showed that the strong type inequality for
$T_W$ follows immediately from Marcinkiewicz interpolation and the
weak type inequality~\eqref{eqn:swm1}. Their estimate is quantitative,
but gives a worse constant in terms of the $\A_p$ characteristic than the best known result.  

\smallskip

Matrix weighted estimates for fractional
integrals have been proved, using the matrix version of the $\A_{p,q}$
weights.  Strong type inequalities were proved by Isralowitz and
Moen~\cite{MR4030471}, who also introduced and proved estimates for a
fractional Christ-Goldberg maximal operator.  Weak type
inequalities were proved in~\cite{MR4859193}.  In both cases the
estimates are quantitative, but it is not known if they are sharp.
Isralowitz and Moen explicitly conjecture that their result is not the best possible.

\section{Convex set-valued functions, matrix weights, and extrapolation}
\label{section:convex}

In this section we consider the problem of extending the Rubio de Francia
extrapolation theorem and the Jones factorization theorem to matrix
weights.  While the proofs of these two results are now relatively
straightforward in the scalar case, there did not seem to be any way
to generalize these proofs to the matrix setting.  In particular, the
fundamental obstacle was generalizing  the
Rubio de Francia iteration algorithm, which required iterates of the
maximal operator.  The Christ-Goldberg operator could not be used, since it maps vector-valued functions to
scalar-valued functions.

The solution was to move beyond vector-valued functions and work with
convex set-valued functions.  The  convex body sparse domination of
Nazarov {\em et al.} was the first indication of the power and utility
of this approach, but in~\cite{mb-2022} the authors developed it considerably.   Before discussing the new techniques required for
the proofs, we first state the two main theorems from~\cite{mb-2022}.

 \begin{theorem}  \label{thm:jones-matrix}
Given $1<p<\infty$ and a matrix weight $W$, $W\in \A_p$
   if and only if there exist commuting matrices  $W_0\in \A_1$ and $W_1 \in \A_\infty$
   such that
   \begin{equation} \label{eqn:jm1}
 W=W_0^{\frac{1}{p}}W_1^{\frac{1}{p'}}.  
\end{equation}
 \end{theorem}   

\begin{remark}
  For simplicity and ease of comparison to the scalar case, we state
  Theorem~\ref{thm:jones-matrix} assuming that the
  matrices $W_0$ and $W_1$ 
  commute.  When we factor an $\A_p$ weight it is always possible to get
  commuting matrices.  For reverse factorization, we can remove this hypothesis, but to do so we must
  replace the product in \eqref{eqn:jm1} with a more complicated
  expression, the geometric  mean of the two matrices:  see~\cite{mb-2022} for details. 
\end{remark}

\smallskip

To state the matrix weight version of Rubio de Francia extrapolation,
we extend our convention about extrapolation pairs to pairs of
vector-valued functions $(\vecf, \vecg) \in \F$ that satisfy
inequalities of the form~\eqref{eqn:extrapol-pairs}.  However,
  we omit the assumption
that they are nonnegative.

\begin{theorem} \label{thm:rubio-matrix}
Given $p_0$, $1\leq p_0\leq \infty$, suppose that for a family of
extrapolation pairs $\F$
and 
for all $W_0\in \A_{p_0}$, the inequality
\[ \|\vecf\|_{L^{p_0}(W_0)} 
  \leq
  N_{p_0}([W_0]_{\A_{p_0}}) \|\vecg\|_{L^{p_0}(W_0)}, \qquad (\vecf,\vecg)\in \F, \]
holds.  Then for all $p$, $1<p<\infty$,  and all $W\in \A_p$, 
\[ \|\vecf\|_{L^{p}(W)} 
  \leq
  N_{p}([W]_{\A_{p}}) \|\vecg\|_{L^{p}(w)}, \qquad (\vecf,\vecg)\in \F, \]
where the function $N_p$ is defined exactly as in Theorem~\ref{thm:rubio}:
\[  N_p ([W]_{\A_{p}}) 
= 
C(n,p_0,p)N_{p_0}\big(c(n,p_0,p) [W]_{\A_{p}} ^{\max\{\frac{p}{p_0},
  \frac{p'}{p_0'}\}}\big).\]
\end{theorem}

\smallskip

\begin{remark}
  As in the scalar case, there are some minor technical difficulties
  to prove weighted norm inequalities for operators using families of
  extrapolation pairs and Theorem~\ref{thm:rubio-matrix}.  For a
  careful discussion, with specific examples, see~\cite{mb-2022}.  
\end{remark}

\smallskip

The proofs of both Theorems ~\ref{thm:jones-matrix}
and~\ref{thm:rubio-matrix} are long and extremely technical.  Our goal
here is to give a conceptual overview of the proofs and the tools
that were developed for them. Underlying the proofs of Theorems~\ref{thm:jones-matrix}
and~\ref{thm:rubio-matrix} was the systematic philosophy of trying to
replicate the proofs of factorization and extrapolation in the scalar
case, especially the elementary proofs of these results
in~\cite{MR2797562,dcu-paseky}.  As we noted above, the proof of
extrapolation in the scalar case depends on four key ideas: the
duality of the $\A_p$ condition, sharp bounds for weighted norm
inequalities for the maximal operator, the Jones factorization
theorem, and the Rubio de Francia iteration operator.  

The duality property of the $\A_p$ condition, as
we noted above, is a consequence of the definition of matrix
$\A_p$ using reducing operators.

The  sharp weighted bounds for the Hardy-Littlewood maximal operator
were more challenging.
A fundamental technical obstacle to the proof was the lack of an
appropriate definition of the maximal operator that can be iterated to produce a
vector-valued function.  To overcome this
problem, the proof passes from vector-valued functions to the larger category
of convex set-valued functions, and uses a convex set-valued
maximal operator. There is actually a
well-developed theory of convex set-valued functions,
see~\cite{MR2458436,CV}.   However,  it does not appear to be
well-known among harmonic analysts.  

We begin with some definitions about convex sets.  For complete
details, see~\cite{mb-2022} and the references it contains.  Let $\K$
denote the family of all convex sets $K\subset \rd$ that are closed,
bounded, and symmetric: i.e., if $x\in K$, then $-x\in K$.  Sometimes
it is necessary also to assume $K$ is absorbing: that is, that $0\in
\mathrm{int}(K)$.  Here, however, we will not worry about this technical
hypothesis.  Given a set $K\subset \rd$, let $|K|=\sup\{ |v| : v \in K
\}$.  Given a matrix $W$, define $WK= \{ Wv : v \in K\}$; note that
$WK$ is also a convex set.   Given two convex sets $K$ and $L$, their
Minkowski sum is the set $K+L = \{ u+ v : u \in K, v\in L \}$, which
is again convex. 

Recall (see~\eqref{eqn:norm-ball}) that  every norm $\rho$ has
associated to it the convex set $K_\rho$
which is its unit ball.  The converse is also true:  to  every convex
set $K$ these is associated to it a unique
norm, $\rho_K$, on $\rd$ such that $K$ is its unit ball.  Moreover, arguing as we did above, there
exists a matrix $W$ such that $\rho_K \approx
\rho_W$.  While generally the proof works with matrices $W\in \A_p$
or with their associated reducing operators,
at a key point it is necessary to work with
more general norms and the underlying convex sets.

A convex set-valued function is a map $F: \rn
\rightarrow \K$.  There are several equivalent ways to define
measurability of such functions $F$; for our purposes a useful and
intuitive definition is that there exists a family
$\{\vecf_k\}_{k=1}^\infty$ of
measurable functions $\vecf_k : \rn \rightarrow\rd$, called selection functions, such that for
almost every $x\in \rn$,
\[ F(x) = \overline{ \{ \vecf_k(x) : k \in \N \} }. \]

Given such a convex set-valued function $F$, we can define for each
$x$ the associated John ellipsoid function $E$.  This ellipsoid-valued
function can be chosen to be measurable. This fact appears to have been known
(see, for instance, Goldberg~\cite{MR2015733}), but a 
proof was not in the literature; a proof was given in~\cite{mb-2022}.

The integral of a convex set-valued function is defined using
selection functions; this is referred to as the Aumann
integral~\cite{MR2458436, MR0185073}.  Given $\Omega
\subset \rn$ and a convex set-valued function $F : \Omega \rightarrow
\K$, define  the collection of integrable selection functions of
$F$,
\[ S^1(\Omega, F) = \{ \vecf \in L^1(\Omega,\rd) : \vecf(x) \in
  F(x)\}.  \]
The Aumann integral of $F$ is defined to be the set
\[ \int_\Omega F(x)\,dx
  =
  \bigg\{ \int_\Omega \vecf(x)\,dx : \vecf \in S^1(\Omega,F)
  \bigg\}. \]
Since $F(x)$ is closed, bounded and convex, the
Aumann integral is also a closed, convex set in $\rd$. (See~\cite{mb-2022}.)

\begin{remark}
  There is a close connection between the Aumann integral and the convex
averages $\langle\langle \vecf \rangle\rangle_Q$ used by Nazarov, {\em
  et al.} to define convex body sparse domination.  Given a
vector-valued function $\vecf$, define $F_{\vecf}(x) =
\clconv\{\vecf(x), -\vecf(x)\}$; that is, the closed convex hull of
the set $\{\vecf(x), -\vecf(x)\}$.  Then $F_{\vecf}$ is a
measurable, convex set-valued function, and
\[ \langle\langle \vecf \rangle\rangle_Q
  = \avgint_Q F_{\vecf}(x)\,dx. \]
Thus, the convex body sparse operator $T_\Ss$ is the infinite
Minkowski sum of  convex set-valued
functions; if $\vecf$ is, for instance, bounded and has compact
support, then this sum converges to a convex set-valued function.
(See~\cite{MR3689742}.)
\end{remark}

\smallskip

We use the Aumann integral to define the convex set-valued maximal
operator.  Given $F : \rn \rightarrow \K$, let
\[ MF(x) = \clconv\bigg(
  \bigcup_Q \avgint_Q F(y)\,dy \cdot \chi_Q (x) \bigg);  \]
that is, $MF(x)$ is the closed, convex hull of the union of the Aumann integral averages of
$F$ over all cubes containing $x$.   Then $MF$ is a measurable,
convex set-valued function.  (See~\cite{mb-2022}.) The intuition
behind this definition is that the Hardy-Littlewood maximal operator
finds the largest average in magnitude, and so uses the supremum.  The
convex set-valued maximal operator finds the largest average in
magnitude in each direction; the union of all these averages preserves
information about both magnitude and direction.

\begin{remark}
  The convex set-valued maximal
operator does not preserve some natural subsets of $\K$.  For
instance, if $F$ is
ellipsoid-valued (that is, $F=W{\mathbf B}$, where $W$ is a matrix
valued function) then $MF$ need not be ellipsoid-valued.
See~\cite{mb-2022} for an example.
\end{remark}

\smallskip

The convex set-valued maximal operator has a number of properties that
correspond to those of the Hardy-Littlewood maximal operator,
replacing inequalities with set inclusion.  For
convex set-valued functions $F$ and $G$, almost every $x\in \rn$, and $\alpha \in
[0,\infty)$, 
\begin{enumerate}

\item $F(x) \subset MF(x)$;

  \item $M(F+G)(x) \subset MF(x) + MG(x)$, where the sum is the
    Minkowski sum;

  \item $M(\alpha F)(x) = \alpha MF(x)$.
\end{enumerate}
The maximal operator $M$ also satisfies $L^p$ norm inequalities.  For
$1\leq p<\infty$, define $L^p_\K(\rn,|\cdot|)$ to be the collection of
all convex set-valued functions $F$ such that
\[ \|F\|_{L^p_\K(\rn,|\cdot|)}
  =
  \bigg(\int_\rn |F(x)|^p\,dx\bigg)^{\frac1p} < \infty.  \]
When $p=\infty$,  define $L^\infty_\K(\rn,|\cdot|)$ by  $\|F\|_{L^\infty_\K(\rn,|\cdot|)}= \||F| \|_{L^\infty(\rn)}<\infty$.
Then for $1<p\leq\infty$, $\|MF\|_{L^p_\K(\rn,|\cdot|)} \leq
C \|F\|_{L^p_\K(\rn,|\cdot|)}$,
and when $p=1$ a weak type inequality holds.

\begin{remark}
  Although the notation $ \|\cdot\|_{L^p_\K(\rn,|\cdot|)}$ is
  deliberately suggestive of a norm, the set $L^p_\K(\rn,|\cdot|)$ is
  not a normed vector space, since the Minkowski sum does not have an
  additive inverse.  It is, however, a complete metric space.
  See~\cite{mb-2022} for details.
\end{remark}

Sharp weighted norm inequalities for the convex maximal operator are
governed by the matrix $\A_p$ weights.  For $1\leq p<\infty$, define
$L^p_\K(\rn,W)=L^p_\K(W)$ to be all convex set-valued functions $F$ such
that
\[ \|F\|_{L^p_\K(\rn,W)}
  =
  \bigg(\int_\rn |W(x)F(x)|^p\,dx\bigg)^{\frac1p} < \infty.  \]
When $p=\infty$, define $L^\infty_\K(\rn,W)$ by
$\|F\|_{L^\infty_\K(\rn,W)}= \|WF\|_{L^\infty_\K(\rn,|\cdot|)}<\infty$.

\begin{theorem} \label{thm:convex-max-wts}
 Given $p$, $1<p\leq\infty$, and $W\in \A_p$, for every $F \in
 L^p_\K(W)$,
 \[ \|MF\|_{L^p_\K(W)}
   \leq
   C(n,d,p)[W]_{\A_p}^{p'}\|F\|_{L^p_\K(W)}. \]
\end{theorem}

The proof of Theorem~\ref{thm:convex-max-wts} uses the $L^p$ norm
inequalities for the Christ-Goldberg maximal operator $M_W$; the quantitative
constant in terms of $[W]_{\A_p}$ was proved by Isralowitz and
Moen~\cite{MR4030471}.  It is the best possible since this is the same
exponent as in the scalar case.  The proof essentially bounds
$\|MF\|_{L^p_\K(W)}$ by the sum of terms of the form
$\|M_W\vecf\|_{L^p(\rn)}$.   In fact there is a converse to this
inequality, showing the close connection between the Christ-Goldberg
maximal operator and the convex set-valued operator:  see
Vuorinen~\cite[Lemma~3.4]{MR4773502}.

\begin{remark}
  An important open problem is to give a direct proof of
  Theorem~\ref{thm:convex-max-wts} that
does not use the bounds for the Christ-Goldberg maximal operator.  In
the scalar case there is a simple proof of the strong $(p,p)$ norm inequalities
for the Hardy-Littlewood maximal operator in Theorem~\ref{thm:scalar-max} due to Christ and
Fefferman~\cite{MR684636} which initially seemed a good model.  However, this proof uses the universal
dyadic maximal operator~\eqref{eqn:universal-max} which does not
appear to generalize to the setting of matrix weights.   We believe
there is a close connection between this problem and the problem of
finding sharp estimates for sparse operators and singular integrals.
\end{remark}

\smallskip

Before we can define the iteration operator, which is used to
construct $\A_1$ weights, 
we need to define the analog of the $\A_1$ condition for convex
set-valued functions.  Given $F : \rn \rightarrow \K$, we say that
$F\in \A_1^\K$ if $MF(x) \subset CF(x)$.  Denote the infimum of all
such constants $C$  by $[F]_{\A_1^\K}$.  There is a very close
relationship between convex set-valued $\A_1^\K$ and matrix $\A_1$.

\begin{theorem} \label{thm:two-A1-defn}
  Given a matrix weight $W$, $W\in \A_1$ if and only if $F=W{\mathbf
    B} \in \A_1^\K$, and $[W]_{\A_1} \approx [F]_{\A_1^\K}$. 
\end{theorem}

The proof of Theorem~\ref{thm:two-A1-defn} in~\cite{mb-2022} required
a careful development of the properties of norm functions and their
duals on the one hand, and the properties of convex sets and their
polar bodies on the other.

With this machinery we can now define the convex set-valued
analog of the Rubio de Francia iteration operator.   Given $F : \rn
\rightarrow \K$  and $W\in \A_p$, $1<p<\infty$, define
\[ \Rdf F(x) = \sum_{k=0}^\infty \frac{M^kF(x)}{2^k
    \|M\|^k_{L^p_\K(W)}}. \]
This sum converges in norm and pointwise in the appropriate
metric.  The iteration operator $\Rdf$ has the following properties which are the exact analogs
of the properties of the scalar iteration operator
given above:
\begin{enumerate}
\item $F(x) \subset \Rdf F(x)$;

\item $\|\Rdf F\|_{L^p_\K(W)} \leq 2 \|F\|_{L^p_\K(W)}$;

  \item $\Rdf F \in \A_1^\K$ and $[\Rdf F]_{\A_1^\K} \leq 2
    \|M\|_{L^p_\K(W)}$.
  \end{enumerate}

  \smallskip

  The final key idea we need to generalize to the matrix setting is the Jones factorization theorem,
  Theorem~\ref{thm:jones-matrix}.    And, in fact, we only need
  reverse factorization, which is very easy to prove in the scalar
  case.  However, in the matrix case, the opposite is true.  Given a
  matrix weight $W \in \A_p$, the construction of $W_0\in \A_1$ and
  $W_1\in \A_\infty$ follows the scalar
proof in~\cite{dcu-paseky}  closely.  There is one  technical
obstacle:  the scalar proof uses a variant of Rubio de Francia
iteration operator based on the 
maximal operator, $M_sf(x) = M(|f|^s)^{\frac1s}$, $s>1$, which is a
sublinear operator.   To define $M_sF(x)$, requires replacing $F$ by an
appropriate ellipsoid (so that the power $F^s$ is well-defined) and then
proving that the resulting operator is sublinear.

The converse,
proving reverse factorization, is very delicate and much more
difficult than in the scalar case.  The proof does not use the definitions of matrix $\A_p$
and $\A_1$ in terms of matrices or reducing operators.
Instead, this proof requires the definition of matrix
$\A_p$ in terms of norms as originally given by Treil and Volberg.
Intuitively, the proof can be thought of as an interpolation argument
between finite dimensional Banach spaces: that is, proving that the
interpolation norm is in fact in $\A_p$.

\smallskip

The proof of extrapolation follows the proof of sharp constant
extrapolation in the scalar case given in~\cite{MR2797562}.  (For a
simpler proof that produces a worse constant that can easily be
adapted to the matrix setting, see also~\cite{dcu-paseky}.) Because of the
way in which the matrix $\A_p$ condition is defined, the proof  also yields extrapolation from the endpoint
$p=\infty$; this gives a quantitative version of an extrapolation result
first proved in the scalar case by Harboure, Mac\'\i as and Segovia~\cite{MR944321}.
(This quantitative version was proved independently 
by Nieraeth~\cite{MR4000248}.)
The proof uses reverse
factorization and the Rubio de Francia iteration operator.  In order
to close the gap between the convex set-valued functions that this
operator creates, and the ellipsoid valued functions that appear in
Theorem~\ref{thm:two-A1-defn} (that is, to pass between matrix $\A_1$
and convex set $\A_1$), the iteration operator is modified so that it
produces ellipsoid-valued functions.    

\begin{remark}
  In the original paper~\cite{mb-2022}, the authors only considered
  vector-valued functions which had real-valued components.  This led
  naturally to considered convex set-valued functions in $\rd$.
  However, in many applications it is important to consider
  vector-valued functions that take values in $\mathbb{C}^d$; to prove
  extrapolation in this setting requires developing a theory of convex
  set-valued functions in $\mathbb{C}^d$.  Kakaroumpas and Soler i
  Gibert~\cite{kakaroumpas2025} have shown that all of the proofs in~\cite{mb-2022} extend to
  this setting; these extensions require a number of 
  technical modifications.  
\end{remark}

\begin{remark}
If we apply Theorem~\ref{thm:rubio-matrix} starting with the sharp
constant for singular integrals when $p=2$, we get that for
$1<p<\infty$ a bound of the form $[W]_{\A_p}^{3\max\{p,p'\}}$.  This is
  always worse, when $p\neq 2$, than the bound obtained
  in~\cite{MR3803292}.  It is unclear why extrapolation produces the
  sharp bounds in the scalar case but not for matrix weights.  
\end{remark}

\smallskip

As in the scalar theory of extrapolation, the matrix theory has been
generalized in several directions.  Here we describe three results.
First, Kakaroumpas and Soler i Gibert~\cite{kakaroumpas2025} have proved that matrix
extrapolation also yields vector inequalities analogous
to~\eqref{eqn:scalar-vv}.  They showed that with the same hypotheses
as Theorem~\ref{thm:rubio-matrix}, you get that
\[  \bigg\| \bigg(\sum_{k=1}^\infty |W\vecf_k|^q
    \bigg)^{\frac{1}{q}} \,dx \bigg\|_{L^p(\rn)}
    \leq
    C  \bigg\| \bigg(\sum_{k=1}^\infty |W\vecg_k|^q
    \bigg)^{\frac{1}{q}} \,dx \bigg\|_{L^p(\rn)},
\]
  when $1< p,\,q<\infty$, $w\in \A_p$, and
  $\{(\vecf_k,\vecg_k)\}_{k=1}^\infty\subset \F$.

  \begin{remark}
    As we noted above, in the scalar case extrapolation also yields classical weak type inequalities.  It is an open problem
    whether extrapolation can be used to prove weak type inequalities
    in the matrix setting.
  \end{remark}

  Second, matrix extrapolation has been extended to weights
  and operators defined with respect to general bases.   Given a basis
  $\B$, we can proceed as we did in the scalar case
  (see~\eqref{eqn:basis-Ap} and\eqref{eqn:basis-max}) and define
  matrix weights $\A_{p,\B}$, $1\leq p\leq \infty$, the
  Christ-Goldberg maximal operator $M_{W,\B}$, and the convex
  set-valued maximal operator $M_\B$.   If we assume $\B$ is a matrix
  Muckenhoupt basis, that is, that $\|M_{W,\B}\vecf\|_{L^p(\rn)} \leq
  C\|\vecf\|_{L^p(\rn)}$ whenever $\W\in \A_p$, then the proof of
  Theorem~\ref{thm:rubio-matrix} goes through with essentially no
  change.  See~\cite{dcu-fs-AFM2025}.

  \begin{remark}
    In~\cite{dcu-fs-AFM2025} the authors also investigate the properties
    of matrix weights.  In the scalar case, while in the definition of
    the $\A_p$ classes we implicitly assumed that $0<w(x)<\infty$
    a.e. (and in fact that $w$ and $w^{-1}$ were  in $L^p_\loc(\rn)$
    and $L^{p'}_\loc(\rn)$), it can be shown that these are necessary
    conditions if the Hardy-Littlewood maximal operator satisfies a
    weak $(p,p)$ inequality.
    (This fact is often glossed over in treatments of the Muckenhoupt
    condition.  See~\cite{garcia-cuerva-rubiodefrancia85} for
    details.)
    By a change of variables, one can redefine the Christ-Goldberg
    maximal operator so that it does not use the inverse $W^{-1}$ in
    its definition.  If one assumes that this operator satisfies a
    weak type inequality, then this implies that $W$ is invertible
    almost everywhere.  See~\cite{dcu-fs-AFM2025} for precise definitions
    and details.
  \end{remark}

  Third, matrix extrapolation can be extended to the off-diagonal
  setting.  The natural conjecture is that extrapolation should
  hold for weights that satisfy a matrix $\A_{p,q}$ condition that
  generalizes the $\A_{p,q}$ weights introduced by Muckenhoupt and
  Wheeden:  for instance, if $1<p\le q<\infty$, $W\in \A_{p,q}$ if
\[
[W]_{\A_{p,q}} = \sup_{Q} \bigg( \avgint_Q \bigg( \avgint_{Q} |W(x)
    W^{-1}(y)|^{p'}_{\op} \,dy
  \bigg)^{\frac{q}{p'}} \,dx \bigg)^{\frac{1}{q}}  < \infty.
\]
However, the proof required a stronger condition.

\begin{theorem} \label{thm:off-diag-extrapol-cubes}
Fix $s$, $1\leq s<\infty$, and a family of extrapolation pairs $\F$.   Suppose that for some $ p_0, q_0, 1 \leq p_0 \leq q_0 \leq \infty $, where $\frac{1}{p_0}-
\frac{1}{q_0}=\frac{1}{s'}$ and $r_0=\frac{q_0}{s}$, and  every matrix weight $W_0$ with $ W_0^s \in \A_{r_0} $,
\begin{equation*} 
    \|  \vecf\|_{L^{q_0}(W_0)}
\leq C([W^s]_{\A_{r_0}})
\| \vecg\|_{L^{q_0}(W_0)}, \qquad (\vecf,\vecg) \in \F.  
\end{equation*}
Then for all $ p,\, q $ such that $ 1 < p \leq q < \infty $, $ \frac{1}{p} - \frac{1}{q} = \frac{1}{p_0} - \frac{1}{q_0}$, $r=\frac{q}{s}$, and for all matrix weights $W$ with $ W^s \in \A_{r} $,
\begin{equation*} 
  \|  \vecf\|_{L^{q}(W)}
\leq  C([W^s]_{\A_{r}})
\| \vecg\|_{L^{q}(W)}, \qquad (\vecf,\vecg) \in \F.  
\end{equation*}
\end{theorem}

The problem that arises is that, unlike in the scalar
case, we no longer have that $W\in \A_{p,q}$ if and only if $W^s\in
\A_r$.  The latter condition is strictly stronger:
see~\cite{dcu-fs-AFM2025} for a counter-example.  And it is this
condition which the proof, generalizing the scalar proof, relies on.
It is an open question whether
Theorem~\ref{thm:off-diag-extrapol-cubes} is true with only the weaker
hypothesis that $W\in \A_{p,q}$.

\begin{remark}
  The proof of Theorem~\ref{thm:off-diag-extrapol-cubes} yields
  quantitative estimates on the constants in terms of $[W^s]_{\A_r}$;
  however, when reduced to the constant exponent case this proof does
  not give the sharp results of Lacey, {\em et al.}~\cite{MR2652182}.
  For a proof which does yield these sharp results, see the
  forthcoming paper~\cite{dcu-fs-2025}.
\end{remark}

\smallskip

The machinery developed to prove Theorems~\ref{thm:jones-matrix}
and~\ref{thm:rubio-matrix} provides additional tools for working with
the convex body sparse domination of Nazarov, {\em et al.}  To
illustrate this,  we conclude this section by
sketching a very recent result in~\cite{DCU-MP-FS2025}.  The authors give a new and
relatively simple proof that convex body sparse operators, and so
singular integrals, are bounded on $L^p(W)$, $1<p<\infty$.    We omit
a number of (minor) technical details to emphasize the overall outline
of the proof, which mimics a proof in the scalar case.  

Fix $p$, a sparse family $\Ss$, and  $\vecf\in L^p(W)$. Then by a duality argument there exists
$\vecg \in L^{p'}(W^{-1})$, $\|\vecg\|_{L^{p'}(W^{-1})} \leq d$,
such that
 \begin{align*}
    \|T_\Ss \vecf\|_{L^p_\K(W)}
   & \leq 2\int_{\rn} \langle T_\Ss\vecf(x), \vecg(x)\rangle \,dx   \\
  & \leq C \sum_{Q\in \Ss}
      \bigg\langle \avgint_Q F_\vecf(y)\,dy, \avgint_Q \vecg(x)\,dx
      \bigg\rangle |E(Q)| \\
   & = C\sum_{Q\in \Ss}
      \int_{E(Q)} \bigg\langle  W(z)\avgint_Q
      F_\vecf(y)\,dy, \\
& \qquad \qquad  \qquad     \avgint_Q
      W^{-1}(z)W(y) W^{-1}(y)\vecg(y)\,dy
                                                                 \bigg\rangle
                                                                 \,dz;
   \\
   \intertext{by the definition of the convex set-valued maximal
   operator and the Christ-Goldberg maximal operator, and by the sharp
   bounds for each,}
    & \leq C\sum_{Q\in \Ss} \int_{E_Q}  |W(z)MF_\vecf(z)|
      M_{W^{-1}}(W^{-1}\vecg)(z)\,dz \\
   & \leq C\|MF_\vecf\|_{L^p_\K (W)}
     \|M_{W^{-1}}(W^{-1}\vecg)\|_{L^{p'}(\rn)}\\
   & \leq C[W]_{A_2}^{p'} \|F_\vecf\|_{L^p_\K(W)}
     [W^{-1}]_{\A_{p'}}^p \|\vecg\|_{L^{p'}(W^{-1})} \\
   & \leq C[W]_{\A_p}^{p+p'} \|\vecf\|_{L^p(W)}.
  \end{align*}

  \begin{remark}
    This proof does not give the best known
    estimates:  when $p=2$ it gives an exponent of $4$ instead of
    $3$.  It is a very interesting question whether this approach can
    be refined to prove the sharp exponent.
  \end{remark}

\section{Further results on matrix weights}
\label{section:beyond}

In this section we briefly summarize some additional results
involving matrix weights that we think are of interest.  We note that
this is a very active field and there are many other results which are
not included here.  Here we restrict to results related to our own work.

We begin with
recent work extending matrix weights to the variable Lebesgue spaces.
The variable Lebesgue spaces $L^\pp$ are a generalization of the
classical $L^p$ spaces, replacing the exponent $p$ by an exponent
function $\pp$.  They have been studied extensively for more than 30
years; for their history and basic
properties, we refer to the
books~\cite{cruz-fiorenza-book,diening-harjulehto-hasto-ruzicka2010}.
Scalar weighted norm inequalities in $L^\pp$ spaces were introduced
in~\cite{MR2927495} (see also~\cite{cruz-diening-hasto2011}), where
bounds for the maximal operator were shown to be controlled by $\A_\pp$
weights, a generalization of the Muckenhoupt condition.  Extrapolation
from $L^{p_0}(w_0)$ to $L^\pp(w)$ was proved in~\cite{MR3572271} (see also~\cite{dcu-paseky}); this
immediately yielded bounds for singular integrals and other operators.

The problem of matrix weights in variable Lebesgue spaces has been
considered in a series of recent papers by Penrod and his
collaborators.  In~\cite{MR4777231} they defined matrix $\A_\pp$
weights and used this to characterize averaging operators.  These they
used to prove that certain convolution operators are bounded and the
related approximate identities converge.  As an application they
generalized the ``$H=W$'' theorem of Meyers and
Serrin~\cite{MR0164252} to matrix weighted space, extending an earlier
result in~\cite{MR3544941}.  In~\cite{dcu-mp-2025} the authors settled
an open problem in the scalar case, proving the existence of a reverse
H\"older inequality for scalar $\A_\pp$ weights.  As an application
they proved that matrix $\A_\pp$ weights satisfy both a left and right
openness property.  This result was new for matrix weights even in the
constant exponent case, and was used by Penrod and
Sweeting~\cite{MP-BS-2026} to prove strong type inequalities for singular integrals.  Very recently, Nieraeth and
Penrod~\cite{zn-mp-2025} used left and right openness to prove that the Christ-Goldberg
maximal operator is bounded on variable Lebesgue spaces, and then used
a version of extrapolation to prove that singular integral operators
are bounded on $L^\pp(W)$ when $W \in \A_\pp$.  A different proof,
that does not use extrapolation, will appear in~\cite{DCU-MP-FS2025}.

Variable Lebesgue spaces are examples of Banach function spaces, and
there has been some interest in extending harmonic analysis and
extrapolation theory to this setting.  In the scalar case,
see~\cite{dcu-paseky,MR2797562}.  Recently, Nieraeth has extended the
techniques in~\cite{mb-2022} to prove a very general extrapolation
theorem in matrix weighted Banach function spaces~\cite{zn-2025}.

\smallskip

In the past 25 years there has been a great deal of work on the
weighted theory for multi\-linear operators.  The weight class
$\A_{\vec{p}}$, $\vec{p}=(p_1,\ldots,p_m)$, was introduced by Lerner,
{\em et~al.}~\cite{MR2483720}; they showed that this class was
necessary and sufficient for the boundedness of a multi-sublinear
maximal operator and nondegenerate multilinear singular integrals.
Their results were extended to the setting of matrix weights by
Kakaroumpas and Nieraeth~\cite{sk-zn2024}.  In the scalar weighted
variable Lebesgue spaces, the results of Lerner, {\em et al.} were
extended to the bilinear setting in~\cite{MR4119260}.  These results
should extend to the multilinear setting, but there are numerous
technical obstacles to overcome in their proof and this remains an
open problem.  Finally, we note that
it is also an open problem to extend the results of Kakaroumpas and Nieraeth 
to variable Lebesgue spaces.

\bibliographystyle{plain}
\bibliography{Convex-Proc-ver-3}

\end{document}